\documentclass[11pt]{article}
\usepackage{latexsym}
\usepackage{amsfonts}
\usepackage{amsmath}
\usepackage{amssymb}

\newtheorem{theorem}{Theorem}
\newtheorem{definition}[theorem]{Definition}
\newtheorem{lemma}[theorem]{Lemma}
\newtheorem{corollary}[theorem]{Corollary}
\newtheorem{proposition}[theorem]{Proposition}
\newtheorem{remark}[theorem]{Remark}
\newtheorem{notation}[theorem]{Notation}
\numberwithin{theorem}{section}
\newcommand{\R}{I\!\!R}
\newcommand{\N}{I\!\!N}

\begin{document}

\title{%%% Canonical regular conditional probability distribution for solutions to martingale problems.\\
A new probabilistic approach to non local and fully non linear  second order partial differential equations}

\author{Jocelyne Bion-Nadal \\ 
CMAP (UMR CNRS 7641)  \\
Ecole Polytechnique  
F-91128 Palaiseau cedex \\
jocelyne.bion-nadal@cmap.polytechnique.fr}

\date{October 5, 2012}

\maketitle

\begin{abstract}
We prove that weakly continuous solutions to martingale problems admit a canonical regular conditional probability distribution.  This allows for the  construction of time consistent convex dynamic procedures in a non dominated setting. Making use of the martingale problem approach for continuous diffusions and diffusions with Levy generator, we give an explicit construction of  such procedures having furthermore a Feller property.\\ 
These procedures   lead to viscosity solution of fully non linear  second order partial differential equations in case of continuous diffusions. In case of diffusions with Levy generator this provides a probabilistic approach for the resolution of non local fully non linear second order PDE.

{\bf MSC: }
%% \end{keywords}
\end{abstract}
%====================================================================================
\section{Introduction} 
\label{sec:intro}

Starting from the solution of the non linear heat equation, 
 Peng \cite{P1,P2} has introduced the notion of $G$-expectation $({\cal E}^G_t)_{t \in \R_+}$ which is a sublinear time consistent dynamic procedure defined on continuous functions on $\Omega$. The state space $\Omega$ is there equal to the set of continuous paths, and ${\cal E}^G_t(X)$ is defined by a stepwise evaluation of the partial differential equation (PDE). Denis at al \cite{DHP} have then proved that ${\cal E}^G_0(X)= \sup_{P^{\theta} \in {\cal P}}E_{P^{\theta}}(X)$, where ${\cal P}$ is a weakly compact set of probability measures which are non dominated (i.e. there is no probability measure $P$ such that $Q \ll P$ for all $Q$ in ${\cal P}$). The properties of $({\cal E}^G_t)_{t \in \R_+}$ are up to a minus sign the properties of a sublinear time consistent dynamic risk measure.
 Since the seminal papers of Peng on $G$-expectation, a new challenge is to develop a theory of time consistent dynamic procedures in a non dominated framework, and also to provide viscosity solutions to fully non linear second order partial differential equations.\\
 Some recent works study the properties of these procedures, in the static case\cite{BNK1} and in the dynamic case \cite{NS,BNK2}.  Other works  construct time consistent dynamic procedures.
The first construction of  $G$-expectation by Peng  \cite{P1,P2} was made on the state of continuous paths from the non linear heat equation 
$-\partial_t v(t,x)-G(D^2v(t,x))=0$, $G(A)=\sup_{\gamma\in {\Theta}}Tr(\gamma \gamma^*A)$ ($\Theta$ bounded) with boundary condition $v(0,x)=f(x)$.  Starting from  the unique viscosity solution of a second order PDE containing a non local term, this construction has been extended in \cite{HP} to produce a time consistent sublinear procedure on the set of c\`adl\`ag paths.  In both cases $G$ is independent of $t$ and $x$ and the $G$-heat equation corresponds to diffusions with no drift and a diffusion coefficient varying between two bounds.  
These works of Peng {\it et al.} have motivated other works
constructing   time consistent dynamic procedures on the set of continuous paths, in a  non dominated framework. 
Nutz has produced time consistent sublinear procedures allowing the $G$ coefficient to vary with both $t$ and $x$ (or equivalently with $\omega$ in the non Markovian case) in \cite{N1}, or allowing both the diffusion coefficient and the drift term to vary between bounds \cite{N2}. In both cases the construction is made $\omega$ by $\omega$. Starting from a given set of probability measures for all $t$ and $\omega$, satisfying some compatibility conditions and analytic properties, Nutz and van Handel \cite{NVH} construct also sublinear time consistent procedures on the space of continuous paths. Soner et al \cite{STZ1,STZ2} have constructed convex (and not only sublinear) time consistent dynamic procedures in the setting of diffusions as well.  Their approach makes use of a solution of a Backward Stochastic Differential Equation (BSDE) associated to every probability measure in the set ${\cal P}$ and then solves an ``aggregation problem''. The diffusion coefficient varies in a domain  independent of $t$ and $x$.  In  \cite{STZ2} it is proved that this construction gives rise to viscosity solutions to fully non linear second order PDE.\\
Time consistent dynamic procedures   on a  filtered probability space or  up to a minus sign  dynamic risk measures have been studied in many papers  \cite{BN03,BN04,CDK,D,KS}. Risk measures are characterized by their dual representation. In particular, in the static case \cite{FS,FR}, $\rho_{0,t}(X) =\sup_{Q \in {\cal P}}(E_Q(-X)-\alpha(Q))$ where ${\cal P}$ is a set of probability measures all absolutely with respect to $P$. Sublinear time consistent  risk measures are fully characterized by a stable set of probability measures \cite{D}. Convex time consistent dynamic risk measures are described by a stable set of probability measures ${\cal P}$ and a local penalty defined on ${\cal P}$ satisfying the cocycle condition \cite{BN03,BN04}.  In the particular case of the Brownian filtration, time consistent dynamic risk measures are limits of solutions of BSDE \cite{DPR}.\\
In order to generalize the construction of time consistent dynamic procedures to the case of non dominated probability measures, the main point is to understand the notion of stable set of probability measures in this new framework. Indeed the usual notion of stable set involves $Q$-conditional expectation.  However $Q$-conditional expectation is defined up to a $Q$-null set  and this is a drawback when one considers non equivalent probability measures. On a Polish space $\Omega$,  for all  sub  $\sigma$- algebra ${\cal B}$ of the Borel $\sigma$ algebra, there exists always a measurable version of the $Q$-condition expectation called regular conditional probability distribution given ${\cal B}$ (\cite{SV4}, Theorem 1.1.6.). However there is no unicity of the $Q$-regular conditional distribution given ${\cal B}$. In order to extend the notion of stability by bifurcation  \cite{BN03} to a   set of non dominated probability measures as well as the local condition for a penalty \cite{BN03}, one needs to be able to make a  coherent  choice for  a measurable version of the conditional expectation.\\   
The first goal of the present paper is to prove that weakly continuous solutions to a martingale problem (Section\ref{secMARPB}) is a very nice setting %%% verifier
in which there exists a  canonical version of the conditional expectation (Section \ref{seccano}). This allows to extend the notion of stability to a  set of probability measures which are non dominated (Section \ref {secstable}). This canonical version has furthermore continuity properties (Section  \ref{secconti}). Specializing to the case of continuous diffusions or to the case of diffusions with Levy generator, we then construct penalties having a Feller property (Section \ref{secTCCC}), satisfying the local condition and the cocycle condition. The canonical version of conditional expectation joint with the Feller penalties  give a generic and constructive method to produce  ${\cal T}$-time consistent convex dynamic procedures on the space of continuous paths and also on the space of c\`adl\`ag paths in a non dominated framework. The set ${\cal T}$ can be any subset of $\R^*_+$ in particular it can be $\R^*_+$ or a discrete subset of  $\R^*_+$. This construction generalizes to the non  dominated framework the construction of \cite{JBN-PDE}. It can be used to construct a great variety of time consistent dynamic procedures in non dominated framework. In the case of continuous paths, the underlying set of probability measures ${\cal P}$ is a ``stable set'' generated by probability measures solution to a martingale problem associated to a continuous diffusion with diffusion coefficient $a(t,x)$ and drift coefficient $b(t,x)$ continuous bounded , $a(t,x)$ invertible for all $(t,x)$ such that $(a,b)$ takes values in a multivalued Borel set. In the case of c\`adl\`ag  paths, the underlying set of probability measures ${\cal P}$ is the ``stable set'' generated by probability measures solutions to a martingale problem associated to diffusions with Levy generator with coefficients  $a(t,x)$ $b(t,x)$ satisfying the above conditions and the jump measure $M(t,x)$ satisfying the hypothesis of \cite{St}, $(a,b,M)$ taking values in a multivalued Borel set. The procedure is defined on the closure (for the norm $\sup_{P \in {\cal P}}E_P(|X|)$ of the lattice vector space of continuous coordinate functions $f(X_{t_1},..X_{t_k})$ ($f$ continuous on $(\R^n)^k$) and also on the cone $f(X_{t_1},..X_{t_k})$, where $f$ is lower semi continuous bounded from below.  A construction for a time consistent procedure on c\`adl\`ag paths is independently proposed in \cite{KPZ}, but within a very 
limited and specific framework: %% the probability measures are also solutions to the martingale problem associated to diffusions with Levy generator. However in \cite{KPZ}, 
the set of probability mesasures is  generated by  probability measures solutions to the martingale problem associated to diffusions with Levy generator where there is no drift, the diffusion coefficient and the jump measure depend only of $t$ (not on $x$),  there is  no multivalued Borel mapping, and no penalty function. Furthermore  the procedure is defined on the uniformly continuous functions on the space of c\`adl\`ag  paths, which does not allow to study PDEs in the context of c\`adl\`ag  paths. \\ 
Our last goal is to prove as in  \cite{JBN-PDE}  that these procedures give a new probabilistic approach to  second order PDE (Section \ref{viscsol}). Indeed, in the case ${\cal T}=\R_+$, these procedures applied to the random variable $h(X_t)$ ($h$ lower semi-continuous bounded from below) give  rise to a time consistent convex Feller process $\tilde h(s,X_s)$. The function $\tilde h$ is lower semi-continuous  on $[0,t] \times \R^n$. Making use of the martingale property of the probability measures, we prove that this lower semi-continuous function is a viscosity supersolution of a   second order PDE: 
\begin{equation}
\left \{\parbox{12cm}{
\begin{eqnarray} 
-\partial_u v(u,x)-f(u, x,Dv(u,x),D^2v(u,x),\tilde Kv(u,x))&=&0 \nonumber \label{eqedp00}\\
v(t,x)& = & h(x)\nonumber 
\end{eqnarray}
}\right.
\end{equation}
The non local term   $\tilde Kv(u,x)(y)=v(u,x+y)-v(u,x)-\frac{y^*\nabla v(u,x)}{1+||y||^2}$ is specific of the case of diffusions with Levy generator. In case of continuous diffusions $f$ depends only on $(u, x,Dv(u,x),D^2v(u,x))$, the PDE is a fully non linear second order partial differential equation. In this last case, 
 our results can be compared to those of \cite{STZ2}. In \cite{STZ2}, the function $f$ can also depend on $v$ but  there are restrictive  conditions on the dual of $f$ which are not needed in our approach. This is due to the fact that the construction in  \cite{STZ2} relies on the existence of solutions to BSDE and also on  the aggregation of these solutions. On the contrary in the present paper we have an explicit construction of the time consistent procedure. Furthermore the stability property of the set of probability measures and the conditions on the penalties imply the time consitency for the procedure (Section \ref{secTCC}).  As mentioned in \cite{JBN-PDE}, the time consistency for the procedure corresponds to the usual dynamic programming principle. Notice also that our probabilistic approach can be considered as a general control problem, where the set of control is a set of non dominated probability measures.  As in \cite{JBN-PDE}, under some conditions the function $\tilde h$ is continuous and is a viscosity solution of the PDE.
% We prove that this procedure applied to the random variable $f(X_t)$ are viscosity supersolution (and under restrictive conditions viscosity solutions) of fully non linear second order PDE.    

%====================================================================================

\section{Martingale problem}
\label{secMARPB}
We fix a finite horizon $T$. In the following for $r \leq T$, the state space 
$\Omega^{r}$   is either the set ${\cal C}([r,T], \R^n)$ of continuous paths on $[r,T]$ endowed with the topology of uniform convergence, or the set ${\cal D}([r,T]$ of c\`adl\`ag paths with the Skorokod topology \cite{Bil}.  Recall that $\Omega^r$ is a Polish space i.e. a complete separable metrizable space \cite{Bil}. \\In both cases $(X_t)_{t \in \R^+}$ denotes the coordinate process, and  ${\cal B}^r_t$ is the $\sigma$-algebra  generated by $\{X_u, \;r \leq u \leq t\} $.
\subsection{Diffusions with continuous coefficients}
\label{secdiff}
As in \cite{JBN-PDE} we consider  the martingale problem introduced in \cite{SV1}. In this subsection the state space is the set of continuous paths, $\Omega^r={\cal C}^r={\cal C}([r,T],\R^n)$. 
Let $a:[0,T] \times \R^n \rightarrow M_n(\R)$ and $b:[0,T] \times \R^n \rightarrow \R^n$ be continuous bounded maps such that for all $(t,x)$,  $a(t,x)$ is  definite positive (i.e. $a$ is strictly elliptic).  For given $\theta \in \R^n$ let 
\begin{equation}
Y^{a,b}_{r,t}(\theta)=exp \{\theta^*(X(t)-X(r))-\int_r^t \theta^*b(u,X(u))  du -\frac{1}{2} \int_s^t \theta^* a(u,X(u)) \theta du \} 
\label{eqmar}
\end{equation}
 Following \cite{SV1} one says that the probability measure $Q^{a,b}_{r,y}$ on $(\Omega^r, {\cal B}^r)$ is  a solution to the martingale problem (\ref{eqmar}) 
starting from $y$ at time $r$ if for all  $\theta \in \R^n$,  $(Y^{a,b}_{r,t}(\theta))_{t \geq r}$ is a martingale on  $(\Omega^r, ({\cal B}^r)_t)_{t \geq r},Q^{a,b}_{r,y})$, and if $Q^{a,b}_{r,y}(\{X_r=y\})=1$.

Denote ${\cal C}^{\infty}_c (\R^n)$ the set of ${\cal C}^{\infty}$ functions with compact support. 
Denote $L^{a,b}_t= \frac{1}{2}\sum_1^{n}a_{ij}(t,x) \frac{\partial^2}{\partial x_i \partial x_j} + \sum_1^{n}b_{i}(t,x)\frac{\partial}{\partial x_i}$. From \cite{SV3} Theorem 2.1, the above martingale problem is equivalent to 
\begin{equation}
Z^{a,b}_{r,t}=f(X_t)-f(X_r)-\int_r^t L^{a,b}_u(f)(X_u)du \;\text{is a $Q^{a,b}_{r,y}$-martingale for all f in} \;{\cal C}^{\infty}_c (\R^n)
\label{eqmart}
\end{equation}
The study of the  martingale problem has been extended in \cite{St} to the case of diffusions with Levy generators. 
\subsection{Diffusions with Levy generators}
\label{secDiLev}
In this subsection  $\Omega^{r}$ is equal to  ${\cal D}^r={\cal D}([r,T],\R^n)$, the set of c\`adl\`ag paths endowed with the Skorokod topology \cite{Bil}. 
\begin{definition}
Hypothesis (M)\\
A Borelian map $M$  defined on $\R_+ \times \R^n$ with values in the set of $\sigma$-finite measures on $\R^n-\{0\}$ satisfies hypothesis (M) if for all Borelian subset $\Delta$ of $\R^n-\{0\}$,
\begin{equation} 
\int _{\Delta}\frac{y}{1+||y||^2}M(s,x,dy)\;\text{is continuous bounded}
\label{eqbound}
\end{equation}
\label{defM}
\end{definition}
Let $K^M_t$ be the operator associated to $M$:
\begin{equation}
K^M_t(f)(x)=\int[f(x+y)-f(x)-\frac{y^*\nabla f(x)}{1+||y||^2}]M(t,x,dy)
\label{eqM}
 \end{equation}
Recall the following result from \cite{St}, Theorem 4.3.

\begin{proposition}
For all $a$ and $b$ continuous bounded and $a$ strictly elliptic, for all $M$ satisfying the hypothesis (M), for all $r$ in $\R_+$ and all $y$ in $\R^n$, there is a unique probability measure $Q^{a,b,M}_{r,y}$ on $({\cal D}^r,{\cal B}^r)$ solution to the martingale problem
\begin{equation}
Z^{a,b,M}_{r,t}=f(X_t)-f(X_r)-\int_r^t (L^{a,b}_u + K^M_u)(f)(X_u)du \;\;\forall f \in {\cal C}^{\infty}_c \R^n)
\label{eqmart2}
\end{equation}
starting from $y$ at time $r$ i. e. such that $Q^{a,b,M}_{r,y}(\{X_r=y\})=1$
\label{propMP2}
\end{proposition}
\subsection{Weakly continuous solution to a martingale problem}
\label{secWCM}
Motivated by the above martingale problem both for continuous diffusions and for diffusions with Levy generators, we consider in the following solutions to a martingale problem in a general setting. A finite horizon $T$ is given. The state space  can be either the set of continuous or of c\`adl\`ag paths. 
%such taht for all $\theta_1, \theta_2 \in \Theta$, for all $r>0$, the map $\theta(u,x)= 1_[0,r[\theta_1 (u,x)+1_[r, \infty[\theta_2(u,x)$ belongs to $\Theta$. 

%We can obviously choose a parameter set $\Delta$ such that Denote $Z^{nu}_{r,t}$

\begin{definition}
 For all   $s \in [0,T[$ and $y \in \R^n$,  let  $Q_{s,y}$ be the unique  solution to the martingale problem $Z$ on $(\Omega^s, {\cal B}^s_T)$, starting from $y$ at time $s$. ($Z$ means a whole family $(Z^i)_{i \in I}$.
\begin{itemize}
\item 
$(Q_{s,y})_{y \in \R^n}$ is  weakly  continuous\\ if $Q_{s,y}$ is a continuous function of $y$ for the weak topology.
\item The martingale problem is additive if for all $i \in I$,
for all $0 \leq r \leq s \leq t$, $Z^i_{r,s}$ is ${\cal B}^r_s$ measurable, $Z^i_{r,s}$ is a right continuous function of $s$, and 
\begin{equation}
 Z^i_{r,t}=Z^i_{r,s}+Z^i_{s,t}\;
\label{eqadd}
\end{equation}
\item The martingale problem is bounded if for all $ 0 \leq r \leq s$, $Z^i_{r,t}$ is bounded.
\end{itemize}
\label{defFellprop}
\end{definition}
%\begin{definition}
%Let $s \geq 0$. Assume that for all   $y \in \R^n$ there is a  solution $Q_{s,y}$ to the martingale problem $(Z_{s,t})_{s \leq t}$ starting from $y$ at time $s$. One says that the solution $(Q_{s,y})_{y \in \R^n}$ to the martingale problem   satisfies  the Markov (resp Feller)  property if for all $f$ continuous  bounded on $(\R^n)^j$, for given $s \leq t_1 \leq..t_k$ the map $y \rightarrow Q_{r,y}(f(X_{t_1},...X_{t_k})$  is a  Borelian (resp continuous) function  on $\R^n$. 

%\label{defFellprop}
%\end{definition}
%\begin{remark}
%The above definition is a definition at given time $r$.\\
 %The Feller property implies obviously the  Markov property.
%\end{remark}
\begin{definition} Hypothesis $H_{\Theta}$

$ \Theta$ is  a  set of Borelian maps 
defined on $\R_+ \times \R^n$ with values in a topological space $E$ such that
\begin{enumerate}
\item for all $\theta \in \Theta$, for all $0 \leq r < t\leq T$ and for all $y$ in $\R^n$, there is a unique probability measure  $Q^{\theta}_{r,y}$ on $(\Omega^r_t,{\cal B}^r_t)$ solution to the additive bounded martingale problem $Z^{\theta}$ starting from $y$ at time $r$.
\item for all $\theta \in \Theta$, and $r <T$, $(Q^{\theta}_{r,y})_{y \in \R^n}$ is  weakly  continuous.
\end{enumerate}
\label{not01}
\end{definition}
Continuous diffusions and diffusions with Levy generators provide examples for $\Theta$:

\begin{proposition}
The notations are those of Section \ref{secdiff}. Let $a$ and $b$ be continuous bounded  on $([0,T] \times \R^n)$ such that  $a(s,x)$ is invertible for all $(s,x)$. Let $Z^{a,b}_{r,t} $ be given by equation (\ref{eqmart}).  The unique solution $Q^{a,b}_{r,y}$ on $({\cal C}^r, {\cal B}^r_T)$  to the additive martingale problem $Z^{a,b}_{r,t}$ starting from $y$ at time $r$ is weakly continuous.\begin{equation}
\Theta=\{(a,b) \;\text{continuous bounded,}\;a \;\text{invertible}\}  
\label{eqThetac}
\end{equation}
satisfies hypothesis $H_{\Theta}$.
\label{propFell}
\end{proposition}
{\bf Proof} The weak continuity follows Theorem 7.1 of  \cite{SV2}). See also Proposition 2.6 of \cite{JBN-PDE} for a detailed argument. It follows easily from the definition that $Z^{a,b}_{r,t}$ is additive and bounded for all $f$ in ${\cal C}^{\infty}_c(\R^n)$.   \hfill $\square $\\\
For diffusions with Levy generators we introduce now another hypothesis.
 \begin{definition}
Hypothesis $M_c$.
The map $M$ with values in $\sigma$ finite measures on $\R^n-\{0\}$ satisfies hypothesis $M_C$ if it satisfies hypothesis $M$ and if 
\begin{equation}
\sup_{s,x}\int [||y||^2 1_{||y|| \leq 1}+||y||1_{||y|| >1}]M(s,x,dy) \leq C
\label{HMc}
\end{equation}
\label{defMc}
\end{definition}
The following Lemma is proved in \cite{LM} (cf the proof of Theorem 20).
\begin{lemma}
Let $A,B,C,D$ be strictly positive.  For all $\epsilon>0$ there is $K>0$ and for all  $\epsilon,\eta>0$, there is $h>0$  such that  for all   $0 \leq r \leq T$,  for all $a$ continuous strictly elliptic bounded by $A$ for all $b$ continuous bounded by $B$ and $M$ satisfying hypothesis $M_C$,  and $||y|| \leq D$,
\begin{equation}
Q^{a,b,M}_{r,y}[\sup_{r \leq v\leq T}||X_v||> K] \leq \epsilon
\label{eqLM0}
\end{equation}
\begin{equation}
\forall u \in [r,T],\;\; Q^{a,b,M}_{r,y}[\sup_{r \leq u \leq v \leq inf(u+h,T)}||X_v-X_u||> \eta] \leq \epsilon
\label{eqLM}
\end{equation}
\label{lemmaLM}
\end{lemma}
Notice that for given $r$ and $y$, the above Lemma implies  the relative weak  compacity of the set of probability measures $Q^{a,b,M}_{r,y}$ satisfying the above conditions. The  equation (\ref{eqLM}) is a ``uniform right convergence''. This property  is stronger than the property  needed to prove the relative weak compacity (Theorem 13.2 of \cite{Bil}). In particular it allows to prove the following proposition.\\ 
If $s_n>s$,  the  probability measure  $Q^{a,b,M}_{s_n,y_n}$ on $\Omega^{s_n}$  can be identified with a probability measure on $\Omega^s$ still denoted $Q^{a,b,M}_{s_n,y_n}$ such that $Q^{a,b,M}_{s_n,y_n}(\{X_u=y_n,\;\;\forall s \leq u \leq s_n\})=1$.

\begin{notation}
Let $0 \leq r \leq s\leq t \leq T$.
\begin{itemize}
\item
Denote  $\pi^r_{[s,t]}$ the canonical projection of $\Omega^r=\Omega^r_T$ onto $\Omega^s_t$.
\begin{equation}
\pi^r_{[s,t]}(\omega)=\omega_{|[s,t]}
\label{eqpirs}
\end{equation}
\item In case  $r=s$ the projection $\pi^r_{[s,t]}$ will be denoted simply $\pi^r_t$
\end{itemize}
\label{nota0}
\end{notation}
\begin{proposition}
Let $a$ be continuous strictly elliptic bounded by $A$, $b$ continuous bounded by $B$ and $M$ satisfying hypothesis $M_C$.  Assume that $(s_n,y_n)$ has the limit $(s,y)$. Then 
\begin{enumerate}
\item 
- If $s_n$ is decreasing to $s$, with the above identification,  $Q^{a,b,M}_{s_n,y_n}$  converges weakly to  $Q^{a,b,M}_{s,y}$. \\
- If $s_n$ is increasing to $s$, the image of   $Q^{a,b,M}_{s_n,y_n}$ by $\pi^{s_n}_{[s,T]}$   converges weakly to $Q^{a,b,M}_{s,y}$. 
\item For all $(t_1,..., t_k)$, $s \leq  t_1<...<t_k  \leq T$, for all $f$ continuous bounded on $(\R^n)^k$, $Q^{a,b,M}_{s_n,y_n}(f(X_{t_1},...,X_{t_k}))$ has the limit $Q^{a,b,M}_{s,y}(f(X_{t_1},...,X_{t_k}))$.
\end{enumerate}
\label{lemmFell}
\end{proposition}
{\bf Proof}

From Lemma \ref{lemmaLM} and from Theorem 13.2 of \cite{Bil}), it follows that\\ $\{Q^{a,b,M}_{s_n,y_n},\; n \in  \N\}$ in case $s_n \geq s$, ($\{(Q^{a,b,M}_{s_n,y_n})(\pi^{s_n}_{[s,T]})^{-1},\; n \in \N\}$) in case $s_n < s$ is relatively compact for the weak topology. The set of probability measures on $(\Omega^s,{\cal B}^s_T)$ is metrizable for the weak topology, thus   there is a subsequence $Q_k$  converging to $Q$ for the weak topology on $\Omega^{s}$. 
As in Section 13 of \cite{Bil} denote $T_Q$ the set of $t \in [s,T]$ for which the projection $\Pi_t:\omega \in \Omega^s \rightarrow \omega(t)$ is continuous except at points forming a $Q$-null set. It is proved in Section 13 of \cite{Bil} that $s,T$ belong to $T_Q$ and that the complement of the set $T_Q$ in $[s,T]$  is at most countable.
\begin{enumerate}
\item Step 1. We prove that equations (\ref{eqLM0}) and (\ref{eqLM}) are also  satisfied by $Q$ for $r=s$. Notice that the weak convergence of $Q_k$ to $Q$  means that  $Q_k(f)$ has the limit $Q(f)$ for continuous functions $f$ but this convergence   is not valid  for general Borelian functions.   Let $T^i_Q$ be an increasing sequence of finite subsets of $T_Q$ containing $s$ and $T$ such that $\cup_i T^i_Q$ is dense in $[s,T]$. It follows from the Mapping Theorem (Theorem 2.7. of \cite{Bil}) and the inequality $R(G) \leq \liminf R_k(G)$ for all open set $G$ and every sequence $R_k$ weakly converging to $R$, that  for all i, 
\begin{equation}
Q[\sup_{s\leq v \leq T,\; v \in T^i_Q}||X_v||> K] \leq \epsilon
\label{eqLM2}
\end{equation}
\begin{equation}
\forall u \in [s,T] \cap T^i_Q,\;\;Q[\sup_{u \leq v \leq (inf(u+h,T),\; v \in T^i_Q}||X_v-X_u||> \eta] \leq \epsilon
\label{eqLM3}
\end{equation}
It follows from the monotone convergence Theorem that one can replace in the above equations $T^i_Q$ by $\cup_i T^i_Q$. Equations (\ref{eqLM0}) and (\ref{eqLM}) follow then for $Q$ making use of the density of $\cup_i T^i_Q$ in $[s,T]$ and of the right continuity of $X_v$ for all $v$. 
\item step 2: We prove that for all $(t_1,..., t_l)$, $s \leq t_1<...<t_l  \leq T$, for all $f$ continuous bounded on $(\R^n)^l$, $Q_k(f(X_{t_1},...,X_{t_l})$ has the limit $Q(f(X_{t_1},...,X_{t_l}))$. \\
From \cite{Bil} Section 13 this is true for all $t_1,..., t_l$, in $T_Q$,  $s \leq t_1<...<t_l \leq  T$.
The set $T_Q$ is dense and contains $s$ and $T$, thus for all $U=(t_1,..,t_l)$  there is a sequence $U^j=(t^j_1,..t^j_l)$ of $l$-uplets $T_Q$ valued decreasing to $U$. The function $f$ being uniformly continuous on compact sets, it follows from equations  (\ref{eqLM0}) and (\ref{eqLM}) satisfied for all $Q_k$ and $Q$, that for all $n$, for all $j\geq J$, all $k$, 
$Q_k(|f(X_{t_1},...,X_{t_l})-f(X_{t^j_1},...,X_{t^j_l})|)<\epsilon$ and 
$Q(|f(X_{t_1},...,X_{t_l})-f(X_{t^j_1},...,X_{t^j_l})|)<\epsilon$. 
From the convergence of $Q_k(f(X_{t^J_1},...,X_{t^J_l})$ to $Q(f(X_{t^J_1},...,X_{t^J_l}))$, we then get the result.
\item Step 3.
Let $s \leq t \leq T$. Let  $Z^{a,b,M}_{s',t}$ be given by equation (\ref{eqmart2}) for some function $f$ ${\cal C}^{\infty}$  with compact support. It follows from hypothesis $M_C$ that the function $(u,x) \in [0,T] \times \R^n \rightarrow L^{a,b}_u(f)(x)+K^M_u(f)(x)$ is continuous bounded. A similar  argument as the above one proves that for all $\epsilon>0$, there are $s \leq s_1 < ..s_p \leq T$ and a continuous bounded function $g_{s,t}$ on $(\R^n)^p$ such that for all $Q$ and $Q_k$, $Q_k(|Z^{a,b,M}_{s,t}-g_{s,t}(X_{s_1}...X_{s_p})|\leq \epsilon$, and $Q(|Z^{a,b,M}_{s,t}-g_{s,t}(X_{s_1}...X_{s_p})|\leq \epsilon$.   
Making use of the martingale property of $(Z^{a,b,M}_{s,t})_{s \leq t \leq T}$ for $Q_k$ it follows that $(Z^{a,b,M}_{s,t})_{s \leq t \leq T}$ is a martingale for $Q$. 
%From the right continuity of $X_s$  it follows that $Z^{a,b,M}_{s,t}$ is a martingale for $Q$.  
\item step 4.
$Q_k$ is weakly converging to $Q$ on $\Omega^{s}$. $X_s$ is continuous on $\Omega^s$, thus for all $\eta>0$, $\{||X_s-y||>\eta\}$ is open. It follows that $Q(\{||X_s-y||>\eta\}) \leq \liminf Q_k(\{||X_s-y||>\eta\})$.\\
-If $s_n<s$,  $Q_k$ is the image of $Q^{a,b,M}_{s_{n_k},y_{n_k}}$. let $\epsilon$, $\eta >0$. Let $h$ such that equation (\ref{eqLM}) is  satisfied for all $Q^{a,b,M}_{s_{n_k},y_{n_k}}$ and $Q$. Let $k_0$ such that for all $k \geq k_0$, $|s_{n_k}-s|<h$ and $||y-y_{n_k}||<\eta$. It follows from (\ref{eqLM}) applied with $Q^{a,b,M}_{s_{n_k},y_{n_k}}$,  $u=s_{n_k}$ and $v=s$ that $Q^{a,b,M}_{s_{n_k},y_{n_k}}(\{||X_s-y||>2\eta\}) < 2 \epsilon$. $Q$ is the weak limit of $Q^{a,b,M}_{s_{n_k},y_{n_k}}(\pi^{s_n}_{[s,T]})^{-1}$, and $s$ belongs to $T_Q$, it follows that $Q(\{X_s=y\})=1$. \\
- If $s_n$ is decreasing to $s$, $Q_k=Q^{a,b,M}_{s_{n_k},y_{n_k}}$. Let $\epsilon>0$, for $k$ large enough, $Q_k(\{|X_s-y|>\epsilon\})=0$.
In both cases this proves that $Q=Q^{a,b,M}_{s,y}$.
\end{enumerate} \hfill $\square $\\
\begin{proposition}
 Let $a$ and $b$ be continuous bounded  on $([0,T] \times \R^n)$ such that  $a(s,x)$ is invertible for all $(s,x)$. Assume that $M$ satisfies the hypothesis $M_C$. $Z^{a,b,M}_{r,t} $  given by equation  (\ref{eqmart2}) is bounded.  The unique solution $Q^{a,b,M}_{r,y}$ on $({\cal D}^r, {\cal B}^r_T)$  to the additive martingale problem $Z^{a,b,M}_{r,t}$ starting from $y$ at time $r$ is weakly continuous.  
The set
\begin{eqnarray}
\Theta=\{(a,b,M) \;a,b\;\text{continuous bounded,}\;a \;\text{invertible},\nonumber\\\;M\;\text{ satisfies hypothesis $M_C$ for some $C>0$} \}  
\label{eqThetad}
\end{eqnarray}
satisfies hypothesis $H_{\Theta}$.

\label{propFellL}
\end{proposition}
{\bf Proof} The existence and unicity follow from \cite{St}. The weak continuity follows then from Proposition \ref{lemmFell} applied with $s_n=r$ for all $n$. As already noticed in step 3 of Proposition \ref{lemmFell}, for all $f \in {\cal C}^{\infty}_c$, and all $s \leq t$, $Z^{a,b,M}_{s,t}$ is bounded. \hfill $\square $.\\ Ths weak continuity for given $r$   can also be found in  \cite{EK2} Lemma 5.2., but there is no detailed proof. \\

%However we give here a different  proof of the continuity. \\
%Let $x_n$ be a sequence with limit $x$ in $\R^n$. The relative weak compacity of the sequence $Q^{a,b,M}_{r,x_n}$ follows from Lemma \ref{LemmaLM}. One can thus consider a subsequence $y_n=x_{\phi(n)}$ such that $Q^{a,b,M}_{r,x_n}$ admits the limt $Q$ for the weak topology.

\section{Canonical representation of the conditional expectation for the  weakly continuous solution to a martingale problem}
\label{seccano}
Let $0 \leq r <s \leq T$. The state space $\Omega^r_T$ is either  ${\cal C}([r,T],\R^n)$ or ${\cal D}([r,T],\R^n)$. The goal of this Section is to prove that in the setting of weakly continuous solutions to martingale problems, there exists a canonical choice for a regular conditional  probability,   and that this choice satisfies  nice  properties.\\
\subsection{Unicity of the solution to a martingale problem for non Markovian parameters}
Let $\Theta$ be a set of parameters satisfying Hypothesis $H_{\Theta}$ (Definition \ref{not01}).
Given $r \in [0,T]$ and $y \in \R^n$, we introduce the following stable set of parameters:
\begin{definition}
Let $0 \leq r \leq T$. ${\Theta}^r_T$ denotes the set of ${\cal B}^r_T$ predictable  processes $\gamma$ on $(\Omega^r,{\cal B}^r_T)$ such that:\\ 
There is a finite subdivision $r=s_0<s_1<...<s_n=T$.\\  For all $i \in \{0,1,...n-1\}$ there is a finite set $I_i$, $I_0=\{1\}$,  a finite partition $(A_{i,j})_{j \in I_i}$ of $\Omega^r$ into  ${\cal B}^r_{s_i}$-measurable sets, and $\theta_0,\;(\theta_{i,j})_{j \in I_i} \in \Theta$ such that 
\begin{equation}
\forall s_i  \leq u  <  s_{i+1}, \;\forall \omega \in \Omega^r,\;\gamma(u,\omega)=\sum_{j \in I_i} \theta_{i,j}(u,X_u(\omega))1_{A_{i,j}}(\omega)
\label{eqstable}
\end{equation}
\label{Slambda}
\end{definition}
We define the associated martingale:
\begin{definition}
For all $\gamma$ in ${\Theta}^r_T$ admitting the above representation, 
denote $(Z^{\gamma}_{r,u})_{r \leq u \leq T}$ the additive process $(Z^{\gamma}_{r,v}=Z^{\gamma}_{r,u}+Z^{\gamma}_{u,v}$ for $r \leq u \leq v \leq T$) defined by 
\begin{equation}
\forall s_i  \leq u  \leq  s_{i+1}, \;\forall \omega \in \Omega,\;Z^{\gamma}_{s_i,u}(\omega)=\sum_{j \in I_i} Z^{\theta_{i,j}}_{s_i,u}(\omega)1_{A_{i,j}}(\omega)
\label{eqstable2}
\end{equation}
\label{defZgamma}
\end{definition}
\begin{lemma}
Assume that  $Q^{\gamma}_{r,x}$ is  a solution to the bounded additive martingale problem $(Z^{\gamma}_{r,t})_{r \leq t \leq T}$ starting from $x$ at time  $r$ (Definition \ref{defFellprop}). Let $(Q^{\gamma}_{r,x})^{s,\omega}$ be a regular conditional distribution of $Q^{\gamma}_{r,x}$ given ${\cal B}^r_s$. There is subset $N$ of $\Omega^r_T$ such that for all $\omega \in \Omega^r_T-N$, $(Q^{\gamma}_{r,x})^{s,\omega}$ is solution on $\Omega^r_T$ to the martingale problem  $(Z^{\gamma}_{s,t})_{s \leq t \leq T}$ and $(Q^{\gamma}_{r,x})^{s,\omega}(\{X_u(\omega')=X_u(\omega),\; \forall r \leq u \leq s\})=1$.
\label{lemmaunic}
\end{lemma}
{\bf Proof} Let $s \leq u \leq v \leq t$.
For all  $B \in  {\cal B}^r_u$ , there is a set $N_{B}$ such that for $Q^{\gamma}_{r,x}(N_{B})=0$ and for all $\omega \in \Omega^r_T - N_{B}$,
\begin{eqnarray}
(Q^{\gamma}_{r,x})^{s,\omega}(Z^{\gamma}_{s,v}-Z^{\gamma}_{s,u})1_B)\nonumber\\=E_{Q^{\gamma}_{r,x}}((Z^{\gamma}_{r,v}-Z^{\gamma}_{r,u})1_B|{\cal B}^r_s)) \nonumber\\
E_{Q^{\gamma}_{r,x}}(E_{Q^{\gamma}_{r,x}}((Z^{\gamma}_{r,v}-Z^{\gamma}_{r,u})1_B|{\cal B}^r_u)|{\cal B}^r_s)=0
\label{eqrc}
\end{eqnarray}
The $\sigma$-algebra ${\cal B}^r_u$  is a Borel $\sigma$-algebra. Thus it is countably generated. Let $A_k$ be a countable family of events  including $\Omega^r_T$ generating ${\cal B}^r_u$.  It follows that there is a subset $N$ of $\Omega^r_T$ such that $Q^{\gamma}_{r,x}(N)=0$ and 
\begin{equation}
(Q^{\gamma}_{r,x})^{s,\omega}(Z^{\gamma}_{s,v}-Z^{\gamma}_{s,u})\xi)=0
\label{eqmartin}
\end{equation}
for all $\xi=1_{A_k}$ and all $\omega \in \Omega^r_T-N$. It follows from the dominated convergence Theorem and the monotone class theorem that equation (\ref{eqmartin}) is satisfied by $1_A$ for all $A$ in ${\cal B}^r_u$. Thus equation (\ref{eqmartin}) is satisfied for all $\xi$ simple ${\cal B}^r_u$ measurable. Another application of the dominated convergence Theorem proves that it is satisfied for all $\xi$  bounded ${\cal B}^r_u$ measurable. We have thus proved that for all $s \leq u \leq v \leq t$, there is a set $N_{u,v}$, $Q^{\gamma}_{r,x}(N_{u,v})=0$ such that  equation (\ref{eqmartin}) is satisfied for all $\xi$  bounded ${\cal B}^r_u$ measurable  and all $\omega \in \Omega^r_T-N_{u,v}$. Let $N=\cup_{u,v \in ({\cal Q}\cap[s,T]) \cup \{s,T\}} N_{u,v}$.  Making use of the rightcontinuity of the process $(Z^{\gamma}_{r,u})_{r \leq u \leq t}$ for given $r$,  equation (\ref{eqmartin}) follows for all $s \leq u \leq v \leq t$,  for all $\xi$  bounded ${\cal B}^r_u$ measurable  and all $\omega \in \Omega^r_T-N$. Thus for all $\omega \in \Omega^r_T-N$(, $Q^{\gamma}_{r,x})^{s,\omega}$ is solution on $\Omega^r_T$ to the martingale problem  $(Z^{\gamma}_{s,t})_{s \leq t \leq T}$. By definition of a regular conditional distribution, $(Q^{\gamma}_{r,x})^{s,\omega}(\{X_u(\omega')=X_u(\omega),\; \forall r \leq u \leq s\})=1$.\hfill $\square $
\begin{proposition} For all process $\gamma$ in $\Theta^r_T$, for all $0 \leq r<t \leq T$ there is at most one probability measure on $(\Omega^r_t,{\cal B}^r_t)$ solution to the additive  martingale problem $(Z^{\gamma}_{r,u})$ starting from $y$ at time $r$.
\label{propuni}
\end{proposition}
{\bf Proof} Let $\gamma$ in $\Theta^r_T$ (Definition \ref{Slambda}). The proof is done by iteration on $n \geq 1$. For $n=1$ it is satisfied by hypothesis on $\Theta$.
Assume that it is satisfied for $n-1$. Let $Q$ be a probability measure solution to the martingale problem $(Z^{\gamma}_{r,u})$ on $(\Omega^r_t,{\cal B}^r_t)$ starting from $y$ at time $r$.  It is easy to see that  the restriction of $Q$ to ${\cal B}^r_{s_{n-1}}$ is solution to the martingale problem $(Z^{\tilde \gamma}_{r,u})_{r \leq u \leq s_{n-1}}$ where $\tilde \gamma$ is the  restriction of $\gamma$ to $[r,s_{n-1}[ \times \Omega^r_t$. By recursion hypothesis, the restriction of $Q$ to ${\cal B}^r_{s_{n-1}}$ is uniquely determined. 
Let $Q^{s_{n-1},\omega}$ be a regular conditional distribution of $Q$ given ${\cal B}^r_{s_{n-1}}$.  From Lemma \ref{lemmaunic}, for all $j \in I_{n-1}$, there is a subset $N_j$ of $A_{n-1,j}$ such that for all $\omega \in A_{n-1,j}-N_j$, $Q^{s_{n-1},\omega}$ is solution to the martingale problem $(Z^{\theta_{n-1,j}}_{s_{n-1},v})_{s_{n-1} \leq v \leq t}$. The function $\theta_{n-1,j}$ belongs to $\Theta$, thus from unicity of the solution to the martingale problem associated to $\theta_{n-1,j}$ starting from $X_{s_{n-1}}(\omega)$ at time $s_{n-1}$, it follows that $Q_{s_{n-1},\omega}$ is uniquely determined for all $\omega \in A_{n-1,j}-N_j$ with $Q(N_j)=0$, for all $j \in I_{n-1}$. Thus $Q$ is uniquely determined.  Therefore there exists at most one probability measure $Q$ solution to the martingale problem $(Z^{\gamma}_{r,u})$ on $(\Omega^r_t,{\cal B}^r_t)$ starting from $y$ at time $r$. \hfill $\square $

\subsection{Construction of a probability measure $\overline Q$ on $(\Omega^r,{\cal B}^r_T)$ associated to a non Markovian parameter}
\label{secovQ}
For every Borelian function $f$ on $\Omega^r_s \times \Omega^s_T$, for all $\omega$ in $\Omega^r_s$, $f^{\omega}$ denotes the Borelian map defined on $\Omega^s_T$ by $f^{\omega}(\omega')=f(\omega,\omega')$.
\begin{lemma}
For  given $s$, let $(Q_{s,y})_{y \in \R^n}$ be  the weakly continuous solution on $(\Omega^s,{\cal B}^s_T$)   to a martingale problem starting from $y$ at time $s$ (Definition \ref{defFellprop}).
Then for all $f$ Borelian bounded  on  $\Omega^r_s \times \Omega^s_T$, the map  \begin{eqnarray}
\Omega^r_s \times \R^n  &\rightarrow&  \R \nonumber\\
(\omega,y) &\rightarrow&Q_{s,y}(f ^{\omega})
\label{eqLC2}
\end{eqnarray}
is measurable when $\Omega^r_s \times \R^n$ is endowed with the  Borel product $\sigma$-algebra.\\
When $f$ is continuous on  $\Omega^r_s \times \Omega^s_T$, the above map is furthermore  separately continuous in each variable.  

\label{lemmaCCR2}
\end{lemma}
{\bf Proof}
\begin{itemize} 
\item Assume that $f$ is continuous bounded.  Being $f^{\omega}$ continuous, the map 
\begin{eqnarray}
y\rightarrow Q_{s,y}(f^{\omega}) \;\;
\label{contFell}
\end{eqnarray}
is continuous by hypothesis.
On the other hand,  $Q_{s,y}$ being a probability measure and  $f$ being continuous bounded, the continuity of the map
\begin{eqnarray}
\omega \rightarrow Q_{s,y}(f^{\omega}).
\label{contFellb}
\end{eqnarray}
  follows from the dominated convergence Theorem.  
Thus the map
\begin{equation}
(\omega,y) \rightarrow Q_{s,y}(f^{\omega})
\label{eqLC2aa}
\end{equation}
is  separately continuous in each variable.    The sets  $\Omega^r_s$ and $\Omega^s_T$ are metrizable and separable, it follows from Lemma  4.50 of \cite{AB} that this map is  jointly measurable for the product Borel $\sigma$-algebra. 
\item Let ${\cal H}$ be the set  of bounded Borelian functions $f$ on $\Omega^r_s \times \Omega^s_T$  such that the map $(\omega,y) \rightarrow Q_{s,y}(f^{\omega})$ is Borelian. ${\cal H}$ is a vector space containing the constant functions.  If $f_n$ is an increasing sequence of non negative functions  in ${\cal C}$ with limit $f$ bounded, it follows from the monotone convergence theorem that $Q_{s,y}(f^{\omega})=\lim_{n \rightarrow \infty} Q_{s,y}(f_n^{\omega})$. It follows that $f$ belongs to ${\cal H}$. From the first step of the proof ${\cal H}$ contains  the class ${\cal C}$ of continuous bounded functions. The class ${\cal C}$ is stable by pointwise multiplication. It follows   from the monotone class theorem as stated in \cite{RY} chapter 0 Theorem 2.2,  that ${\cal H}$ contains all the bounded Borelian functions on 
$\Omega^r_s \times \Omega^s_T$. \hfill $\square $
\end{itemize}

\begin{corollary}
For all $f$ Borelian bounded  on  $\Omega^r_s \times \Omega^s_T$, the map  
\begin{equation}
\Omega^r_s  \rightarrow  \R 
\label{eq}
\end{equation}

\begin{equation}
\omega \rightarrow Q_{s,X_s(\omega)}(f^{\omega})
\label{eqLC2c}
\end{equation}
is Borelian.
\label{corbor}
\end{corollary}
{\bf Proof} The map $\omega \in \Omega^r_s \rightarrow X_s(\omega) \in \R^n$  is Borelian, thus the  corollary follows from Lemma \ref{lemmaCCR2}. \hfill $\square $
\begin{proposition}
 Let $0  \leq r <s$. Let $J$ be a finite subset of $\N$. For given $j \in J$, let $(Q^{j}_{s,y})_{y \in \R^n}$ be a weakly continuous family of probability measures on $\Omega^s_T$. Let $(h_j)_{j \in J}$ be non negative Borelian  maps on $\R^n$ such that $\forall y \in \R^n,\;\; \sum_{j \in J}h_j=1$. 
 Let $Q_1$ be a probability measure on $(\Omega^{r}_s, {\cal B}^{r}_s)$.\\
 There is a unique probability measure $\tilde Q$ on $(\Omega^r_s \times \Omega^s_T ,{\cal B}^r_s \times {\cal B}^s_T)$  such that  for all $f$ Borelian bounded on $\Omega^r_s \times \Omega^s_T $,
\begin{equation}
 \tilde Q(f)=\sum_{j \in J}\int_{\Omega^{r}_s}h_j(X_s(\omega)) Q^j_{s,X_s(\omega)}(f^{\omega})dQ_1(\omega)
\label{eqrc00}
\end{equation}
\label{propproba}
\end{proposition}
{\bf Proof} 
For all $f$ Borelian bounded it follows from Corollary  \ref{corbor}  that the 
map 
\begin{equation}
T(f):\omega \in \Omega^{r}_s \rightarrow \sum_{j \in J} h_j(X_s(\omega)) Q^j_{s,X_s(\omega)}(f^{\omega})
\label{eqT}
\end{equation}
 is Borelian. It is also bounded. Thus $\int_{\Omega^{r}_s} \sum_{j \in J}h_j(X_s(\omega)) Q^j_{s,X_s(\omega)}(f^{\omega})dQ_1(\omega)$ is well defined. The map 
\begin{equation}
L:\; f \rightarrow \sum_{j \in J}\int_{\Omega^{r}_s}h_j(X_s(\omega)) Q^j_{s,X_s(\omega)}(f^{\omega})dQ_1(\omega)
\label{eqrc01}
\end{equation}
defines  a non negative linear form on the vector space of continuous bounded functions on  $\Omega^r_s \times \Omega^s_T$. 
Let $f_n$ be a decreasing  sequence of continuous bounded functions on $\Omega^r_s \times \Omega^s_T$ converging to $0$. From the dominated convergence Theorem, it follows that for all $j$ and $\omega$,  $Q^j_{s, X_s(\omega)}(f_n^{\omega})$ is a decreasing bounded sequence with limit $0$. Applying again the dominated convergence Theorem, it follows that $L(f_n)$ is decreasing  with limit $0$. 
From Daniell Stone Theorem it follows that there is a unique probability measure $\tilde Q$ on $(\Omega^r_s \times \Omega^s_T ,{\cal B}^r_s \times {\cal B}^s_T)$,   such that  for all $f$ continuous  bounded, equation (\ref{eqrc00}) is satisfied. 
The equality (\ref{eqrc00}) follows then for all Borelian bounded function as in the end of the proof of Lemma \ref{lemmaCCR2}. \hfill $\square $

\begin{definition} Denote $\overline Q$ the probability measure on $\Omega^r_T$ deduced from $\tilde Q$ by the Borelian map 
\begin{eqnarray}
i:\Omega^r_s \times \Omega^s_T& \rightarrow &\Omega^r_T \nonumber\\
 (\omega, \omega')& \rightarrow & i(\omega, \omega')
\label{eqtansfer1}
\end{eqnarray}
where
\begin{eqnarray}
i(\omega, \omega')(u)&=&\omega(u),\; \text{for}\; r \leq u \leq s\nonumber\\
&=&\omega'(u)+\omega(s)-\omega'(s),\;\text{for}\; s \leq u \leq T
\label{eqtansfer2}
\end{eqnarray}
\label{deftransfer}
\end{definition}
\subsection{Canonical regular conditional probability for $\overline Q$ given ${\cal B}^r_s$}

Recall that $\Pi^r_s$ denotes the canonical projection from  $\Omega^r_T$ to $\Omega^r_s$. (cf Notation \ref{eqpirs}).

\begin{proposition}
\begin{itemize}
\item 
 For all $\psi$ Borelian bounded on $\Omega ^r_T$, for all $\omega$ in $\Omega^r_T$, let 
\begin{equation}
T_s(\psi)(\omega)=\sum_{j \in J}h_j(X_s(\omega))(Q^j)_{s,X_s(\omega)}((\psi \circ i)^{\pi^r_s(\omega)})
\label{eqTQ}
\end{equation}
The map $T_s(\psi)$ defined on $\Omega ^r_T$ is ${\cal B}^r_s$ measurable.
\item  For all $\tilde Q$, and all $\psi$ Borelian bounded,  $T_s(\psi)$ is a 
${\cal B}^r_s$ measurable version of the $\overline Q$ conditional expectation of $\psi$ given ${\cal B}^r_s$. \\ 
For all $g$  ${\cal B}^r_s$ measurable, $T_s(g)=g$.
\end{itemize}
\label{propL2} 
\end{proposition}
{\bf Proof}
\begin{itemize}
\item
The map  $i$ is  measurable for the Borel $\sigma$-algebras on both side, it follows from Lemma \ref{lemmaCCR2} that for every Borelian map $\psi$  on $\Omega^r_T$, the map 
\begin{eqnarray}
\Omega^r_s \times \R^n  &\rightarrow&  \R \nonumber\\
(\omega,y) &\rightarrow&Q^j_{s,y}((\psi \circ i)^{\omega})
\label{eqLC2b00}
\end{eqnarray}
is Borelian for all $j$. Notice that the map $\pi^r_s$ is measurable from $\Omega ^r_T$ endowed with the $\sigma$-algebra ${\cal B}^r_s$ into $\Omega ^r_s$ endowed with its Borel $\sigma$-algebra. Composing the above map (\ref{eqLC2b00}) with 
\begin{eqnarray}
\Omega^r_T \times \R^n  &\rightarrow& \Omega^r_s \times \R^n  \nonumber\\
(\omega,y) &\rightarrow& (\pi^r_s (\omega),y)
\label{eqLC2bb}
\end{eqnarray}
it follows that 
\begin{eqnarray}
\Omega^r_T \times \R^n  &\rightarrow&  \R \nonumber\\
(\omega,y) &\rightarrow&Q^j_{s,y}((\psi \circ i)^{\pi^r_s(\omega)})
\label{eqLC2bbb}
\end{eqnarray}
 is measurable for the $\sigma$-algebra product of ${\cal B}^r_s$ and of ${\cal B}(\R^n)$. $h_j$ a Borelian map on $\R^n$ and $X_s$ is ${\cal B}^r_s$ measurable,  it follows that $T_s(\psi)$ is ${\cal B}^r_s$ measurable for all Borelian map $\psi$.
\item 
Notice that every  ${\cal B}^r_s$ measurable function $g$ defined on  $\Omega^r_T$   can be factorized  $g=\tilde g\circ\pi^r_s$ for some $\tilde g$ ${\cal B}^r_s$ measurable defined on $\Omega^r_s$. This  result is deduced from the monotone class theorem (\cite{RY}) applied with the class ${\cal C}$ of functions $f(X_{t_1},X_{t_2},..,X_{t_k})$ ,$k \in \N^*$, $f$ Borelian on $(\R^n)^k$, $r \leq t_1 \leq  t_2 \leq t_k \leq s$.  For all $\omega$ in $\Omega^r_T$ and $\omega'$ in $\Omega^s_T$, $\pi^r_s \circ i(\pi^r_s(\omega),\omega')=\pi^r_s(\omega)$. It follows that for all $g$ ${\cal B}^r_s$ measurable, the map $(g\circ i)^{\pi^r_s(\omega)}$ is constant equal to $g(\omega)$. From the equality  $\sum_{j \in J}h_j=1$, it follows that for all $g$ ${\cal B}^r_s$ measurable, $T_s(g)=g$.\\
$T_s(\psi)$ factorizes,  $T_s(\psi)=\tilde T_s(\psi) \circ \pi^r_s$ where $\tilde T_s(\psi)$ is defined on $\Omega^r_s$ by 
\begin{equation}
\tilde T_s(\psi)(\omega)=
\sum_{j \in J}h_j(X_s(\omega)) Q^j_{s,X_s(\omega)}(\psi \circ i)^{\omega}
\label{eq002}
\end{equation}
From the  definition of $\overline Q$ (Proposition \ref{propproba} and Definition \ref{deftransfer}), it follows that  for all $\Psi$ ${\cal B}^r_T$ measurable  , and $g=\tilde g \circ \pi^r_s$  ${\cal B}^r_s$ measurable defined on $\Omega^r_T$, 
\begin{equation}
 \overline Q(\psi g)=\sum_{j \in J}\int_{\Omega^{r}_s}(\tilde g)(\omega)h_j(X_s(\omega)) Q^j_{s,X_s(\omega)}(\psi \circ i)^{\omega})dQ_1(\omega)
\label{eqrc001}
\end{equation}
From the equality $T_s(\psi)=\tilde T_s(\psi)\circ \Pi^r_s$, it follows that 
\begin{equation}
 \overline Q(T_s(\psi) g)=\int_{\Omega^{r}_s}(\tilde g)(\omega)(\tilde T_s(\psi)(\omega)dQ_1(\omega)
\label{eqrc003}
\end{equation}
The map  $T_s(\psi)$ being ${\cal B}^r_s$ measurable, it follows from equations (\ref{eqrc001}), (\ref{eq002}) and (\ref{eqrc003}) that $T_s(\psi)$ is a ${\cal B}^r_s$ measurable version of the $\overline Q$ conditional expectation of $\psi$ given ${\cal B}^r_s$. \hfill $\square $
\end{itemize}
\subsection{Existence of a solution to the martingale problem for non Markovian parameters}
The goal of this Section is to prove that for all $\gamma \in \Theta^r_T$ (Definition \ref{Slambda}), the martingale problem $Z^{\gamma}$ introduced in Definition \ref{defZgamma} has a solution.
\begin{proposition} 
The hypothesis are the same as in Proposition \ref{propproba}.  
Assume furthermore that $Q^{1}=Q^{1}_{r,x}$ is  the unique solution  on $\Omega^r_s$ to  the $(Z^{\theta}_t)_{r\leq t}$ martingale problem starting from $x$ at time $r$. Assume that for all $j \in J$,  $Q^{j}_{s,y}$ is  the unique solution  on $\Omega^{s}_T$ to  the $(Z^{\theta_j}_t)_{s \leq t}$  martingale problem starting from $y$ at time $s$. Assume that there is a finite partition $(A_j)_{j \in J}$ of $\Omega^r$ in ${\cal B}^r_s$ measurable sets such that for all $j \in J$, $h_j=1_{A_j}$. Assume that  for all $j$ in $J$, $(Q^{j}_{s,y})_{y \in \R^n}$ is weakly continuous.  The probability measure $\overline Q$ is  the unique  solution to the $Z^{\beta}$ martingale problem starting from $x$ at time $r$ with $\beta(u,x)=\theta(u,x),\;$ for all $r \leq u <s$ and   $\beta(u,x)=\sum_{j \in J}1_{A_j}\theta_j(u,x)$,\; for all $s \leq u$. 
\label{propcondprob}
\end{proposition}
{\bf Proof}

\begin{itemize}
\item Let $s \leq t<u \leq T$. By definition of $\beta$ and the additive property of the martingales, $Z^{\beta}_{r,u}-Z^{\beta}_{r,t}=Z^{\beta}_{t,u}= \sum_{j \in J}1_{A_j}
(Z^{\theta_j}_{s,u}-Z^{\theta_j}_{s,t})$. For all $\omega$ in $\Omega^r_s$ and $\omega'$ in $\Omega^s_T$, such that $X_s(\omega)=X_s(\omega')$, for all $s \leq v \leq T$, $i(\omega,\omega')(v)=\omega'(v)$, and $[(Z^{\theta_j}_{s,u}-Z^{\theta_j}_{s,t}) \circ i]^{\omega}(\omega')= (Z^{\theta_j}_{s,u}-Z^{\theta_j}_{s,t})(\omega')$. 
Notice also that for all $\psi$ ${\cal B}^r_t$ measurable and all $\omega$ in $\Omega^r_s$, the map $(\psi \circ i)^{\omega}$ is ${\cal B}^s_t$ measurable. Thus for all $\psi$ ${\cal B}^r_t$ measurable 
\begin{eqnarray}
1_{A_j}(X_s(\omega))Q^j_{s,X_s(\omega)}[([\psi (Z^{\beta}_{r,u}-Z^{\beta}_{r,t}]\circ i)^{\omega}]\nonumber\\
=1_{A_j}(X_s(\omega))Q^j_{s,X_s(\omega)} [1_{\{\omega'|\; X_s(\omega)=X_s(\omega')\}}(\psi \circ i)^{\omega}[(Z^{\theta_j}_{s,u}-Z^{\theta_j}_{s,t}) \circ i]^{\omega}]\nonumber\\
=1_{A_j}(X_s(\omega)Q^j_{s,X_s(\omega)}[(\psi \circ i)^{\omega}(Z^{\theta_j}_{s,u}-Z^{\theta_j}_{s,t})]
\label{eqmart1}
\end{eqnarray}
$(Z^{\theta_j}_{s,t})_{s \leq t \leq T}$ is a $Q^j_{s,y}$ martingale for all $y$, it follows that the last term of (\ref{eqmart1}) is equal to $0$. It follows then from (\ref{eqmart1}) and the definition of $\overline Q$ (Section \ref{secovQ}) that $Z^{\beta}_{r,t}$ is a $\overline Q$ martingale for $s \leq t \leq T$.
\item Let $r \leq t<u \leq s$. $Z^{\beta}_{r,u}-Z^{\beta}_{r,t}=(Z^{\theta}_{r,u}-Z^{\theta}_{r,t})=(Z^{\theta}_{r,u}-Z^{\theta}_{r,t})\circ (\pi^r_s)$. Let $\psi$ be a function ${\cal B}^r_t$ measurable on $\Omega^r_T$, $\psi=\tilde \psi \circ (\pi^r_s)$, where $\tilde \psi$ is a function  ${\cal B}^r_t$ measurable defined on $\Omega^r_s$. Being 
$\psi$ and $Z^{\beta}_{r,u}-Z^{\beta}_{r,t}$ ${\cal B}^r_s$ measurable, it follows from the second part of the proof of Proposition \ref{propL2} that 
$\overline Q(\psi (Z^{\beta}_{r,u}-Z^{\beta}_{r,t}))=Q_1(\tilde \psi (Z^{\theta}_{r,u}-Z^{\theta}_{r,t}))$. Being $(Z^{\theta}_{r,t})_{r \leq t \leq s}$ a $Q_1$ martingale, it follows that $Z^{\beta}_{r,t}$ is a $\overline Q$ martingale for $r \leq t \leq s$. \hfill $\square $
\end{itemize}
\subsection{Canonical regular conditional probability}
\begin{notation} 
For all $r \leq s$, for all $\omega \in \Omega^r_T$ and $\omega' \in \Omega^s_t$ such that $X_s(\omega)=X_s(\omega')$ we denote $\omega*\omega'$  the element of $\Omega^r_T$ such that 
\begin{eqnarray}
\omega*\omega'(u)&=&\omega(u)\;\; \forall r \leq u \leq s\nonumber\\
\omega*\omega'(u)&=&\omega'(u)\;\; \forall s <u \leq T
\label{eq*}
\end{eqnarray}
\label{not*}
\end{notation}

\begin{proposition}
\begin{enumerate}
\item
Let  $r \geq 0$ and $y$ in $\R^n$.  Let $\theta \in \Theta$. For all $s \in [r,T[$ and  all bounded  map $\psi$ ${\cal B}^r_T$-measurable defined on $\Omega^r$ ,  the map $T^{\theta}_s(\psi)$ defined on $\Omega^r$ by $T^{\theta}_s(\psi)(\omega)=Q^{\theta}_{s, X_s(\omega)}((\psi \circ i)^{\pi^r_s(\omega)})$ is  a ${\cal B}^r_s$-measurable   version of the $Q^{\theta}_{r,y}$-conditional expectation of $\psi$ given ${\cal B}^r_s$.
We call it the canonical  regular conditional distribution of $Q^{\theta}_{r,y}$.
\item 
\begin{eqnarray}
T^{\theta}_s(\psi)(\omega)&=&Q^{\theta}_{s, X_s(\omega)}((\psi \circ i)^{\pi^r_s(\omega)})\label{eq39a}\\&=&\int_{\Omega^s_T} \psi(\omega*\omega')1_{X_s(\omega)=X_s(\omega')}dQ^{\theta}_{s, X_s(\omega)}(\omega')
\label{eqTs}
\end{eqnarray}
\item It satisfies the following chain rule for all map $\phi$ ${\cal B}^{r_0}_T$-measurable defined on $\Omega^{r_0}$,
\begin{equation}
\forall r_0\leq  r \leq  s \leq T,\;\; T^{\theta}_r(\phi)=T^{\theta}_r(T^{\theta}_s(\phi))
\label{eqcomp}
\end{equation}
\end{enumerate}
\label{proCRC}
\end{proposition}
{\bf Proof}
\begin{enumerate}
\item
Apply Proposition \ref{propcondprob} with  $Q_1=Q^{\theta}_{r,y}$, $J=\{j\}$, and $Q^j_{s,x}= Q^{\theta}_{s,x}$. It follows that the probability measure $\overline Q$ satisfies the $Z^{\theta}$ martingale problem  on $\Omega^r_T$ starting from $y$ at time $r$. By unicity of the solution to the martingale problem it follows that $\overline Q=Q^{\theta}_{r,y}$.  The result follows then from Proposition \ref {propL2}.
\item It follows immediately from the definitions of $i$, $\Pi^r_s$ and $*$ that 
\begin{equation}
1_{\{X_s(\omega)=X_s(\omega')\}}(\psi \circ i)^{\pi^r_s(\omega)}(\omega')=1_{\{X_s(\omega)=X_s(\omega')\}}\psi(\omega*\omega')
\label{eq*2}
\end{equation}
The equation (\ref{eqTs}) follows then from $Q^{\theta}_{s, X_s(\omega)}(1_{\{X_s(\omega)=X_s(\omega')\}})=1$
\item Let $\overline Q$ be as in 1.. From  the definition of $\overline Q$, Definition \ref{deftransfer}, equations (\ref{eqrc00}) and (\ref{eqTs}) it follows from the equality   $\overline Q=Q^{\theta}_{r,y}$ that for all function $\psi$ ${\cal B}^{r}_T$-measurable defined on  $\Omega^{r}$ ,
\begin{equation}
\int_{\Omega^r_T} \psi (\omega')dQ^{\theta}_{r,y}(\omega')=\int_{\Omega^r_s}[ \int_{\Omega^s_T}\psi(\omega_1*\omega_2) 1_{X_s(\omega_1)=X_s(\omega_2)}dQ^{\theta}_{s, X_s(\omega_1)}(\omega_2)]
dQ^{\theta}_{r,y}(\omega_1)
\label{eqcr1}
\end{equation}
Let $\phi$ ${\cal B}^{r_0}_T$ measurable defined on $\Omega^{r_0}$. Let $\omega \in \Omega^{r_0}_r$. 
Applying the above equation to the function $\psi$  defined on $\Omega^r$ by $\psi(\omega')=\phi(\omega*\omega')1_{\{X_r(\omega)=X_r(\omega')\}}$, it follows from equation (\ref{eqTs})  that 
\begin{eqnarray}
T^{\theta}_r(\phi)=\int_{\Omega^r_T} \phi (\omega*\omega')1_{X_r(\omega)=X_r(\omega')}dQ^{\theta}_{r,X_r(\omega)}(\omega')=\nonumber\\
\int_{\Omega^r_T}1_{X_r(\omega)=X_r(\omega_1*\omega_2
)}dQ^{\theta}_{r,X_r(\omega)}(\omega_1) [ \int_{\Omega^s_T}\phi(\omega*\omega_1*\omega_2)1_{X_s(\omega_1)=X_s(\omega_2)} dQ^{\theta}_{s, X_s(\omega_1)}(\omega_2)]\nonumber\\
=\int_{\Omega^r_T}1_{X_r(\omega)=X_r(\omega_1)T^{\theta}_s(\phi)(\omega*\omega_1)
)}dQ^{\theta}_{r,X_r(\omega)}(\omega_1) 
\label{eqcr2}
\end{eqnarray}
This proves the chain rule equation (\ref{eqcomp}). \hfill $\square $
\end{enumerate}

\begin{theorem}
Let $0 \leq r \leq T$. For all $\gamma \in \Theta^r_T$ and $y$ in $\R^n$, there is a unique solution to the martingale problem $Z^{\gamma}$ on $(\Omega^r,{\cal B}^r_T)$ starting from $y$ at time $r$. We denote it $Q^{\gamma}_{r,y}$. Let $r\leq s <T$. Consider the expression of $\gamma$ as in  Definition \ref{Slambda}.  For all $r \leq s<T$, let $i$ such that $s_i \leq s <s_{i+1}$.  For all $\psi$ defined on $\Omega^r$ bounded  ${\cal B}^r_{s_{i+1}}$ measurable,  the map  $T^{\gamma}_s(\psi)$ defined as 
\begin{equation}
T^{\gamma}_s(\psi)=\sum_{j \in J_i}1_{A_{i,j}}T^{\theta_{i,j}}_s(\psi)
\label{Tgamma}
\end{equation}
is $ {\cal B}^r_{s_i}$ measurable. 
The map   $T^{\gamma}_s$ admits a unique extension to all bounded ${\cal B}^r_T$ measurable functions such that the chain rule is satisfied: 
\begin{equation}
\forall r \leq s \leq t \leq T,\;  T^{\gamma}_s(\psi)=T^{\gamma}_s (T^{\gamma}_t(\psi))
\label{eqcomp3a}
\end{equation} 
$T^{\gamma}_s(\psi)$ is a  ${\cal B}^r_s$ measurable version of the $Q^{\gamma}_{r,y}$ conditional expectation of $\psi$ given ${\cal B}^r_s$. It is called the canonical regular conditional distribution of $Q^{\gamma}_{r,y}$ given ${\cal B}^r_s$.  
\label{propcanreg}
\end{theorem}
{\bf Proof} 
For given $r$ and $y$ the existence of $Q^{\gamma}_{r,y}$ follows by induction from  Proposition \ref{propcondprob}. The unicity follows from  Proposition \ref{propuni}. Let $s_i$ be as subdivision asociated to $\gamma$ as in Definition \ref{Slambda}.
Let $ s_i \leq  s <s_{i+1}\leq T$. Consider the probability measures $\overline Q$ constructed on $\Omega^r_{s_{i+1}}$  from $Q^{\gamma}_{r,y}$ on $\Omega^r_s$ and from $\sum 1_{A_{i,j}}(\omega)Q^{\theta_{i,j}}_{s,X_s(\omega)}$ on $\Omega^s_T$ as in Section \ref{secovQ}. From Proposition \ref{propcondprob} it follows that $\overline Q$ solves the martingale problem for $\gamma$ on $\Omega^r_{s_{i+1}}$ starting from $y$ at time $r$. From the unicity it follows that the restriction of $Q^{\gamma}_{r,y}$ to ${\cal B}^r_{s_{i+1}}$ is equal to $\overline Q$. From Proposition \ref{propL2} equation (\ref {eqTQ}) and equation (\ref{eq39a}), it follows that for all  $\psi$  ${\cal B}^r_{s_{i+1}}$ measurable, 
$T^{\gamma}_s(\psi)$ defined by equation (\ref{Tgamma}) is a ${\cal B}^r_s$ measurable  version of the $\overline Q$- conditional expectation of $\psi$ given ${\cal B}^r_s$. It is thus also a ${\cal B}^r_s$ measurable  version of the $Q^{\gamma}_{r,y}$-conditional expectation of $\psi$ given ${\cal B}^r_s$.\\
Making use of the chain rule for $T^{\theta}$, Proposition \ref{proCRC}, it follows that $T^{\gamma}_s$ admits a unique extension to all ${\cal B}^r_T$ measurable functions such that the chain rule (\ref{eqcomp3a}) is satisfied.  It is given  by a stepwise evaluation: for $s_i \leq s < s_{i+1}$, $T^{\gamma}_s(\psi)=T^{\gamma}_s (T^{\gamma}_{s_{i+1}}(...T^{\gamma}_{s_{n-1}}(\psi))$. The result follows then by induction.
\hfill $\square $
\begin{corollary}
\begin{enumerate}
\item
Let $r\geq 0$, for all $\gamma \in \Theta^r_T$, and $r \leq s <T$, the non negative linear map $T^{\gamma}_s$ is continuous from below.  
It can be uniquely  extended to the Borelian functions  bounded from below by the formula 
\begin{equation}
T^{\gamma}_s(\psi)=lim_{n \rightarrow \infty}T^{\gamma}_s(\psi \wedge n)
\label{eqext}
\end{equation}
This extension is continuous from below.
\item Equation (\ref{eqTs}) (resp (\ref{Tgamma})) is still satisfied for random variables bounded from below, for all $\theta \in \Theta$ (resp $\gamma \in \Theta^r_T$).
\item 
$T^{\gamma}_s(\psi)$ is  a ${\cal B}^r_s$-measurable   version of the $Q^{\gamma}_{r,y}$-conditional expectation of $\psi$ given ${\cal B}^r_s$.
\item
The chain rule is satisfied for all function $\psi$ ${\cal B}^r_T$-measurable  bounded from below, 
\begin{equation}
\forall  r \leq  s \leq t \leq T,\;\; T^{\gamma}_s(\psi)=T^{\gamma}_s(T^{\gamma}_t(\psi))
\label{eqcompe}
\end{equation}
\end{enumerate}
\label{corext}
\end{corollary}
{\bf Proof} The continuity from below for $T^{\theta}_s$, $\theta \in \Theta$ follows  from equation (\ref{eqTs}) and from the  monotone convergence theorem applied to the probability measure $Q^{\theta}_{s, X_s(\omega)}$ for given $\omega$. It follows for general $T^{\gamma}_s$ from the formula (\ref{Tgamma}). The other statements follow from the continuity from below of $T^{\gamma}_s$ and   from the monotone convergence theorem for probability measures (or for conditional expectations to prove {\it 3.}).  \hfill $\square $

We are now able to give a definition of a stable set of probability measures for probability measures which are non dominated by some probability measure.

\section{Stable set of probability measures and penalties}
\label{secstable}
We concentrate now on the non dominated framework. As already noticed the usual definition of a  stable set of probability measures cannot be used because the $Q$-conditional expectation is defined up to $Q$-null set.  The existence of a canonical choice for the regular conditional probability is crucial in the construction of dynamic procedures in non dominated setting.
\subsection{Stable set of probability measures}
\label{susecsta} 
Let $0 \leq r \leq T$. Let $\Omega^r$ be ${\cal C}^r$ or ${\cal D}^r$ with the canonical filtration. For all $r \leq s\leq T$, denote $({\cal B}_b)^r_s$ the set of bounded functions defined on $\Omega^r$, ${\cal B}^r_s$-measurable.  Let ${\cal T}$ be a subset of $[r,T]$.
\begin{definition}
A set ${\cal Q}$ of probability measures on $(\Omega^r,{\cal B}^r_T)$ is  ${\cal T}$ stable if 
\begin{enumerate}
\item Special choice of a regular conditional distribution\\
For all $Q \in {\cal Q}$ for all $r\leq  s \leq T$  $s \in  {\cal T}$,  there is a  non negative linear map continuous from below 
$T^Q_s: ({\cal B}_b)^r_T \rightarrow ({\cal B}_b)^r_s$  such that for all $\psi
\in ({\cal B}_b)^r_T$, $T^Q_s(\psi)$ is a version of the $Q$-conditional expectation of $\psi$ given $({\cal B}_b)^r_s$. 
\item Chain rule\\
For all map $\psi$  in $({\cal B}_b)^r_T$,  for all $s,t$ in ${\cal T}$, $s \leq t$ 
\begin{equation}
\;\; T^{Q}_s(\psi)=T^{Q}_s(T^{Q}_t(\psi))
\label{eqcomp3}
\end{equation}
\item Stability by composition\\
For all $Q$ and $R$ in ${\cal Q}$, for all $s \in  {\cal T}$  there is a probability measure $S$ in ${\cal Q}$ such that  
\begin{eqnarray}
\forall s \leq u \leq T, \;
T^S_u= T^R_u \nonumber\\ 
\forall r \leq u<s, \forall \phi \in ({\cal B}_b)^r_s,\; T^S_u(\phi)= T^Q_u(\phi)
\label{eqcompo} 
\end{eqnarray}
\item stability by bifurcation:\\
for all $ s \in  {\cal T}$ for all $Q$ and $R$ in ${\cal Q}$,  for all $A \subset \Omega^r$,   $A \in {\cal B}^r_s$, there is a probability measure $S$ in ${\cal Q}$ such that 
\begin{equation}
\forall s \leq u \leq T, \;T^{S}_u=1_AT^{Q}_u+1_{A^c}T^{R}_u
\label{eqbif}
\end{equation}
\end{enumerate}
 \label{defstable}
\end{definition}
In all the following to simplify the notations we assume that ${\cal T}=[r,T]$.
\begin{theorem}
Assume that  $\Theta$ satisfies hypothesis $H_{\Theta}$ (Definition \ref{not01}). 

Let $0 \leq r \leq T$ and $y$ in $\R^n$.  Let $\Theta^r_T$ be the set introduced in Definition \ref{Slambda}. The set 
\begin{equation}
{\cal Q}^{\Theta}_{r,y}=\{Q^{\gamma}_{r,y},\gamma \in  \Theta^r_T\}
\label{eqsta}
\end{equation}
is a stable set of probability measures. \\
More precisely, let $Q=Q^{\gamma}_{r,y}$ and $R=Q^{\delta}_{r,y}$  in ${\cal Q}^{\Theta}_{r,y}$.
\begin*
\item{i)}
 Let $\lambda(u,\omega)=\gamma(u,\omega)$ for all $r \leq u <s$ and $\lambda(u,\omega)=\delta(u,\omega)$ for all $s \leq u < T $. The probability measure $S=Q^{\lambda}_{r,y}$  satisfies (\ref{eqcompo}).
\item{ii)}
Let $\eta(u,\omega)=\gamma(u,\omega)$ for all $r \leq u < s$ and $\eta(u,\omega)=1_A(\omega)\gamma(u,\omega)+1_{A^c}\delta(u,\omega)$ for all $s \leq u < T $. The probability measure $S=Q^{\eta}_{r,y}$  satisfies (\ref{eqbif}).
\end*
\label{thmsta}
\end{theorem}
{\bf Proof}
For $\gamma \in \Theta^r_T$ and $r \leq s < T$, let $T^{Q^{\gamma}_{r,y}}_s=T^{\gamma}_s$. It follows  from Corollary \ref{corext} that $T^{Q^{\gamma}_{r,y}}_s$ is linear non negative and continuous from below. Properties {\it 1.} and {\it 2.} follow from Theorem \ref{propcanreg}. \\
- Proof of property {\it 3.}
 The process $\lambda$ defined in ${i)}$ belongs to $\Theta^r_T$ and from Theorem  \ref{propcanreg} there is a unique  probability measure $Q^{\lambda}_{r,y}$ solution on $(\Omega^r,{\cal B}^r_T)$ to the martingale problem associated to $\lambda$ starting from $y$ at time $r$.   It follows also from equation (\ref{Tgamma}) and the chain rule that  for all $r \leq u  \leq v \leq T$ for all $\psi$ defined on $\Omega^r$ ${\cal B}^r_v$ measurable, the ${\cal B}^r_u$  measurable map  $T^{\lambda}_u(\psi)$ depends only on the restriction of $\lambda$ to $[u,v[$. This end the proof of  the stability by composition.\\
- Property {\it 4.} is proved in the same way considering the process $\eta$ defined in ${ii)}$\hfill $\square $\\
\linebreak
As in \cite{JBN-PDE} we can furthermore consider a multivalued Borel mapping $\Gamma$ from $\R^+ \times\R^n$ into $E$ (definition 2.8 of \cite{JBN-PDE}). For all $(t,x)$, $\Gamma(t,x)$ is a subset of $E$ containing $0$. Let $\Theta(\Gamma)=\{\theta \in \Theta\;|\theta(t,x) \in \Gamma(t,x) \;\forall (t,x)\}$.   We then get the following corollary
\begin{corollary}
 Assume that  $\Theta$ satisfies hypothesis $H_{\Theta}$.  For all multivalued Borel mapping $\Gamma$, the set 
\begin{equation}
{\cal Q}^{\Theta}(\Gamma)_{r,y}=\{Q^{\gamma}_{r,y},\gamma \in  \Theta(\Gamma)^r_T\}
\label{eqstab}
\end{equation}
is a stable set of probability measures. We call it the stable set of probability measures generated by $\Theta (\Gamma)$. 
\label{corsta}
\end{corollary}

\subsection{penalties}
\label{subsecpen}
We want to construct  time consistent dynamic convex Feller procedures, (and not only sublinear). Therefore as for the construction of time consistent dynamic risk measures on a filtered probability space, we need to introduce penalties satisfying some properties (\cite{BN03}).  To encompass the non dominated framework, we must adapt the definition of the required properties for penalties.
\begin{definition} 
 Assume that  $\Theta$ satisfies hypothesis $H_{\Theta}$. Let $\Gamma$ be a  multivalued Borel mapping.  Let ${\cal Q}={\cal Q}^{\Theta}(\Gamma)_{r,y}$.  A penalty defined on ${\cal Q}$ is a family $(\alpha_{s,t}(Q))_{r \leq s \leq t\leq T}$ of ${\cal B}^r_s$-measurable  functions on $\Omega^r$ ${\cal B}^r_s$-measurable bounded from above such that 
\begin{enumerate}
\item It is local:\\
For all $\gamma,\eta \in {\Theta(\Gamma)^r_T}$, for all $A \in {\cal B}^r_s$, if $1_A \gamma(u,\omega)=1_A \eta(u,\omega)$, for all $u \in[s,t[$, then $1_A \alpha_{s,t}(Q^{\gamma}_{r,y})=1_A \alpha_{s,t}(Q^{\eta}_{r,y})$. 
\item It satisfies the cocycle condition:
\begin{equation}
\forall Q \in {\cal Q},\;\forall r \leq s \leq t \leq u \leq T,\;
\alpha_{s,u}(Q)=\alpha_{s,t}(Q)- T^Q_s(-\alpha_{t,u}(Q))
\label{eqcc}
\end{equation}
\end{enumerate}
\label{defpen}
\end{definition}
 \begin{remark}
Penalties could be defined for more general stable sets  of probability measures. 
\label{rempen}
\end{remark}
The construction that we have made in \cite{JBN-PDE} in the setting  of (equivalent) probability measures  solution to a martingale problem associated to a continuous diffusion leads also to penalties  in the non dominated framework.  This will be detailed when we construct time consistent dynamic procedures (Section \ref{secTCCC}).\\
Recall also that one of our goals is to construct time consistent convex procedures which give rise to viscosity solutions to second order partial differential equations.  Therefore as in \cite{JBN-PDE} we are interested in continuity  properties.  
\section{Continuity properties for the canonical conditional probability}
\label{secconti}
\subsection{A continuity property in the case of continuous paths}
\label{secccp}

\begin{proposition}
 Case of ${\cal C}$:
 For all $u$ in $\R^+$ and $x \in \R^n$, let $Q_{u,x}$ be the unique solution to a martingale problem. Assume that for some $s$, $(Q_{s,y})_{y \in \R^n}$ is weakly continuous. 
\begin{enumerate}
\item 
Then for all $0 \leq r <s$  for all $f$  bounded uniformly continuous (for the uniform norm) on $\Omega^r_T$, the canonical conditional probability  $T^{Q_{r,x}}_s(f):$
\begin{eqnarray}
\Omega^r_T   &\rightarrow&  \R \nonumber\\
(\omega &\rightarrow&Q_{s,X_s(\omega)}((f \circ i)^{\pi^r_s(\omega}))=Q_{s,X_s(\omega)}((f(\omega*\omega')1_{X_s(\omega)=X_s(\omega')})
\label{eqCCCa}
\end{eqnarray}
is a continuous version of the $Q_{r,x}$ conditional expectation given ${\cal B}^r_s$. 
\item
If $f=g \circ \pi^r_{t'}$ for some $g$ continuous on $\Omega^{r}_{t'}$ with compact support. Let $t''=inf(s,t')$. Then $Q_{s,X_s(\omega)}((f \circ i)^{\pi^r_s(\omega}))=h \circ \pi^r_{t"}$ for some continuous function $h$ on $\Omega^{r}_{t"}$ with compact support.
\end{enumerate}
\label{propCCCa}
\end{proposition}
{\bf Proof}
\begin{enumerate}
\item 
The map $i$ introduced in Definition \ref{deftransfer} is uniformly continuous for the uniform norm, it follows that for all $f$ uniformly continuous bounded on $\Omega^r_T$, the map $f \circ i$ is uniformly continuous bounded on the product space  
$\Omega^r_s \times \Omega^s_T$ for the uniform norm. 
Being $f \circ i$ uniformly continuous, for all $\epsilon>0$ there is $\eta>0$ such that for $||\omega-\omega'||<\eta$, $||(f \circ i)^{\omega}-(f \circ i)^{\omega'}||< \epsilon$. The continuity of the map 
 \begin{eqnarray}
\Omega^r_s \times \R^n  &\rightarrow&  \R \nonumber\\
(\omega,y) &\rightarrow&Q_{s,y}((f \circ i)^{\omega})
\label{eqLC2b}
\end{eqnarray}
follows then easily from the weak continuity of the family $(Q_{s,y})_{y \in  \R^n}$.
being $X_s$ continuous from $\Omega^r_T$ to $\R^n$ and $\pi^r_s$ continuous from $\Omega^r_T$ to $\Omega^r_s$, the continuity follows by composition. From Proposition \ref{proCRC}, equation (\ref{eqCCCa}) describes the canonical regular conditional distribution.\\
\item Assume now that $f=g \circ \pi^r_{t'}$ for some $g$ continuous on $\Omega^{r}_{t'}$  with compact support $K$. 
\begin{enumerate}
\item If $t'\leq s$, $f$ is  ${\cal B}^r_s$ measurable, and from Proposition \ref{propL2}, $ T^{Q_{r,x}}_s(f)=f$
\item If $s \leq t'$, then $t''=s$. The maps $g$ and $\pi^r_{t'}$ are uniformly continuous thus from 1., $T^{Q_{r,x}}_s(f)$ is continuous and ${\cal B}^r_s$ measurable. Therefore it can be written $T^{Q_{r,x}}_s(f)=h\circ \pi^r_s$ for some continuous function $h$.
The map $\pi:\omega \in \Omega^{r}_{t'} \rightarrow \omega_{|[r,s]} \in \Omega^{r}_{s}$ is continuous. Thus $\pi(K)$ is a compact subspace of $\Omega^{r}_{s}$.  If $\omega \notin  \pi(K)$ and $\omega' \in \Omega^s_T$, $Q_{s,X_s(\omega)}((f(\omega * \omega')1_{ X_s(\omega)= X_s(\omega')})=0$.  This proves the result with the support of $h$ contained in  $\pi(K)$. 
\end{enumerate}
\end{enumerate}
\hfill $\square $

Notice that this provides a class of continuous functions  stable by $T^{\theta}$:
\begin{notation}
In case $\Omega^r={\cal C}^r$, we denote 
\begin{equation}
{\cal V}^r_t=\{f\circ \pi^r_u,r \leq u \leq t, f \text{continuous on $\Omega^r_u$ with compact support }\}
\label{eqV}
\end{equation}
For all $Y \in {\cal V}^r_t$, for all $y \in \R^k$ and $s \geq r$,
$T^{Q_{r,y}}_s(Y)$ belongs to ${\cal V}^r_t$
\label{notcor}
\end{notation}

\begin{lemma}
Let $\Omega$ be a Polish space. Let ${\cal Q}$ be a set of probability measures on $(\Omega,{\cal B}(\Omega))$.  Assume that ${\cal Q}$ is tight. For all $f$ continuous bounded on ${\cal C}^r$, there is an increasing sequence of continuous functions with compact support $f_n$ such that $f=\lim_{n \rightarrow \infty}f_n \; Q \; a.s. \; \forall Q \in {\cal Q}$.
\label{lemmatig} 
\end{lemma}
{\bf Proof}
By hypothesis there is an increasing sequence ${\cal K}_n$ of compact subspaces of $\Omega$ such that $Q({\cal K}_n^c)< \frac{1}{n}$ for all $Q \in {\cal Q}$. Let $\phi_n$ be an increasing sequence of continuous functions with compact support, $1_{{\cal K}_n} \leq \phi_n \leq 1$.
 Assume that $f$ is continuous bounded, let $g=\lim_{n \rightarrow \infty} f\phi_n$. $0 \leq  ||f-f\phi_n||\leq||f||$.  $Q(|f-f\phi_n|)\leq||f||\frac{1}{n}$. From the dominated convergence Theorem, it follows that   $f=g=\lim f \phi_n\;\; Q \; a.s. \; \forall Q \in {\cal Q}$.

\subsection{Feller property}
\label{secfell}

\begin{definition} 
\begin{enumerate}
\item
Let $s$ in $\R^+$. For all  $x$ in $\R^n$, let $Q^{\theta}_{s,x}$ be the unique weakly continuous solution to a martingale problem starting from $x$ at time $s$. \\
$(Q^{\theta}_{s,x})_{x \in \R^n}$ has  the Feller property if for all $t\geq s$, for all $f$ continuous bounded on $\R^n$,  the map $x \rightarrow Q^{\theta}_{s,x}(f(X_t))$ is continuous. We denote it $T_{st}(f)$. 
\item
Assume that $\Theta$ satisfies hypothesis $H_{\theta}$. $\Theta$ has the Feller property if  every $(Q^{\theta}_{s,x})$ has  the Feller property.
\end{enumerate}
\label{defFeller}
\end{definition}
The sets $\Theta$ constructed in section 2 have the Feller property:
\begin{lemma}
The sets  $\Theta$  given by equation (\ref{eqThetac}) in case of continuous diffusions and by equation (\ref{eqThetad}) in case of diffusions with Levy generator satisfy the property $H_{\theta}$ and have the Feller property.
\label{lemFellprop}
\end{lemma} 
{\bf Proof} The property $H_{\theta}$ has been  proved in  Proposition \ref{propFell} in case of continuous diffusions. The map $X_t$ is continuous, thus the Feller property is satisfied. In case of diffusions with Levy generator, the property $H_{\theta}$ has been  proved in  Proposition \ref{propFellL}. The Feller property has been proved in  Proposition \ref{lemmFell}.
\begin{lemma}
Case of ${\cal D}$ or ${\cal C}$. let $Q^{\theta}_{s,x}$ be the unique weakly continuous solution  to a martingale problem starting from $x$ at time $s$. 
Assume that  $(Q^{\theta}_{s,x})_{x \in \R^n}$ has the Feller property.  Then 
\begin{equation}
\omega \in \Omega^r_T \rightarrow T_{st}(f)(X_s)(\omega)
\label{eqDCC2}
\end{equation}
 is the canonical  version of the $Q^{\theta}_{r,x}$ conditional expectation of $f(X_t)$ given ${\cal B}^r_{s}$. 
\label{lemmaCC2}
\end{lemma}
{\bf Proof}
Let $\psi=f(X_t)$, $t \geq s$. With the notations of Proposition \ref{proCRC}, $1_{X_s(\omega)=X_s(\omega')}(\psi \circ i)^{\pi^r_s(\omega)}(\omega')=f(X_t)(\omega')$. The Lemma follows then easily from Proposition \ref{proCRC}, and the definition of $T_{st}(f)$.
\begin{notation}
Let ${\cal Q}$ be a weakly relatively compact set of probability measures. Denote as in \cite{DHP} and \cite{BNK1} $L^1(c)$ the Banach space obtained from the completion and separation of continuous function for the c norm $c(f)=\sup_{Q \in {\cal Q}}Q(|f|)$.
\label{notc}
\end{notation} 

\begin{lemma}
Let ${\cal P}$ be a set of probability measures on ${\cal C}([s,T],\R^n)$. Assume that ${\cal P}$ is tight. Let ${\cal T}$ be a dense subset of $[s,T]$. The set of continuous bounded coordinates map $f(X_{t_1},... X_{t_k})$, $s\leq t_1<...<t_k\leq T$, $t_i \in {\cal T}$ 
and $f \in {\cal C}_b((\R^n)^k)$ is dense in $L^1(c)$.

\label{lemmacont}
\end{lemma}

{\bf Proof}
Being ${\cal P}$ tight,  there is a compact subset $K$ of $\Omega^s={\cal C}([s,T], \R^n)$ such that for all $P \in {\cal P}$, $P(\omega^s-K)< \epsilon$. Let 
${\cal E}=\{f(X_{t_1},... X_{t_k}),\;s\leq t_1<...<t_k\leq t,\;t_i \in {\cal T},\; f\in {\cal C}_b((\R^n)^k)\}$. The restriction to $K$ of elements of ${\cal E}$ is an algebra which separates the points of the Haussdorf  compact space $K$. Thus from Stone Weierstrass theorem it is dense in ${\cal C}(K)$. This proves the result.  
%Let $\Omega^s_t= {\cal C}_b([s,t], \R^n)$

\begin{notation}
Let $0 \leq r \leq s\leq T$. For all $j \in \N^*$ and $k \in \N^*$, for all $f$ Borelian on $(\R^n)^{j+k}$, given $r \leq s_1<...<s_j<s \leq t_1<..t_k \leq T$, denote $f_{(s_1,..,s_j),(t_1,...t_k)}$ the map defined on $\Omega^r_T$ by $f_{(s_1,..,s_j),(t_1,...t_k)}(\omega)= f(X_{s_1},..X_{s_j},X_{t_1},...X_{t_k})(\omega)$.
\label{notf}
\end{notation}

\begin{proposition}
Case of ${\cal D}$ or ${\cal C}$.
For all $s$ in $\R^+$ and $x \in \R^n$, let $Q^{\theta}_{s,x}$ be the unique weakly continuous solution to a martingale problem starting from $x$ at time $s$. Assume that for all $s$ $(Q^{\theta}_{s,x})_{x \in \R^n}$ has the Feller property.
Then for all $f$ continuous bounded (resp. Borelian bounded) on $(\R^n)^{j+k}$, there is a continuous bounded (resp. Borelian bounded) map $\hat f$ on $(\R^n)^{j+1}$ such that for all $r <s$ and $x$ in $\R^n$,
\begin{equation}
\omega \in \Omega^r_T \rightarrow \hat f((X_{s_1}(\omega),...,X_{s_j}(\omega),X_s(\omega))
\label{eqDCC2b}
\end{equation}
 is the canonical  version of the $Q_{r,x}$ conditional expectation of $f_{(s_1,..,s_j),(t_1,...t_k)}$ given ${\cal B}^r_{s}$, i.e. $T^{\theta}_s(f_{(s_1,..,s_j),(t_1,...t_k)})=\hat f_{s_1,...,s_j,s}$
\label{propDCC2}
\end{proposition}

{\bf Proof} 
We prove the proposition   for all $f$  continuous bounded on $(\R^n)^{j+k}$. The result follows then for $f$ Borelian bounded from the monotone class theorem.
\begin{itemize}
\item
Assume that $k=1$. 
%The hypothesis means that for given $t>s$ for all  $h$ continuous bounded function on $\R^n$ , the map 
%\begin{equation}
%y \rightarrow Q_{s,y}(h(X_t))
%\label{eqTs1}
%\end{equation}
%is continuous.
Let $\psi=f_{(s_1,..,s_j),(t_1)}$. With the notations of Prop  \ref{proCRC}, $1_{\{ \omega'|\;\omega'(s)=\omega(s)\}}(\psi \circ i)^{\pi^r_s(\omega)}(\omega')=fX_{s_1}(\omega),...,X_{s_j}(\omega),X_{t_1}(\omega'))$. 
\begin{eqnarray}
Q_{s,X_s(\omega)}((\psi \circ i)^{\pi^r_s(\omega)})=Q_{s,X_s(\omega)}(1_{\{\omega'|\; \omega'(s)=\omega(s)\}}(\psi \circ i)^{\pi^r_s(\omega)})\nonumber\\
=Q_{s,X_s(\omega)}(1_{\{\omega'|\; \omega'(s)=\omega(s)\}}(\phi_{\omega}(X_{t_1}))
\label{eqscont}
\end{eqnarray}
where $\phi_{\omega}$ is the continuous function $\phi_{\omega}(z)=f(X_{s_1}(\omega),...,X_{s_j}(\omega),z)$.   It follows from the Feller property  that for given $x_1,...,x_j$, the map 
\begin{equation}
y \rightarrow Q_{s,y}(f(x_1,...x_j,X_{t_1})\;\text{ is continuous}
\label{eqyconti}
\end{equation} 
 We prove now  that  
the map $\hat f : (x_1,...x_j,y) \rightarrow Q_{s,y}(f(x_1,...x_j,X_{t_1}))$ is continuous.\\The family of probability measures $(Q_{s,y})_{y \in \R^n}$ on $\Omega^s_T$ is weakly continuous thus for given $K>0$, $(Q_{s,y})_{||y|| \leq K}$ is tight. Thus for all $\epsilon>0$, there is a compact subset ${\cal K}$ of $\Omega^s_T$ such that for $||y|| \leq K$, $Q_{s,y}({\cal K}^c) < \epsilon$. $\cal K$ is  compact  for the Skorohod topology, thus from Theorem 12.3 of \cite{Bil}  there is $A>0$ such that $\sup_{\omega' \in {\cal K}}|X_{t_1}(\omega')| \leq A$.
Let $K'>0$. Due to the uniform continuity of $f$ on the compact set $\{||x_i||\leq K'\;\forall 1 \leq i \leq j,\; ||x_{t_1}|| \leq A\}$, being $Q_{s,y}$ a probability measure, for all $\epsilon>0$ there is $\eta>0$ such that for $||x_i|| \leq K'$, $||x'_i|| \leq K'$ and  $||x_i-x'_i||<\eta$,  
\begin{equation}|Q_{s,y}(f(x_1,..,x_j,X_{t_1}))-Q_{s,y}(f(x'_1,..,x'_j,X_{t_1}))|<2(1+||f||) \epsilon, \;\;\forall y, \; \;||y|| \leq K
\label{equnifc}
\end{equation} 
 The continuity of $\hat f:\;(x_1,...x_k,y) \rightarrow Q_{s,y}(g(x_1,...x_j,X_{t_1}))$  follows then easily from (\ref{eqyconti}) and (\ref {equnifc}). 
\item It  follows  from equation (\ref{eqscont}) and Proposition \ref{proCRC} that the map $\omega \rightarrow \hat f((X_{s_1}(\omega),...,X_{s_j}(\omega),X_s(\omega))$ 
is the canonical  version of the $Q_{r,x}$ conditional expectation of $f_{(s_1,..,s_j),(t_1,...t_k)}$ given ${\cal B}^r_{s}$.
\item By iteration. Assume that the result is satisfied for $k-1 \geq 1$. We prove it for $k$. It follows from the first step applied with $s=t_{k-1}$ that there 
is a continuous function $\psi$ on $\R^{j+(k-1)}$ such that $\psi_{(s_1,..,s_j),(t_1,...t_{k-1})}(\omega)$  is the canonical  version of the $Q_{r,x}$ conditional expectation of $f_{(s_1,..,s_j),(t_1,...t_k)}$ given ${\cal B}^r_{t_{k-1}}$,  which means that \\$T^{\theta}_{t_{k-1}}(f_{(s_1,..,s_j),(t_1,...t_k)}=\psi_{(s_1,..,s_j),(t_1,...t_{k-1})}$.  The result follows then by iteration making use of the composition rule  (\ref{eqcomp}).
\end{itemize} 
\hfill $\square $
\begin{notation}
Given $r\leq s$, denote 
${\cal W}^r_s$ the algebra of continuous coordinates maps ${\cal B}^r_s$ measurable.
\begin{equation}
{\cal W}^r_s=\{f_{s_1,..,s_j}|r \leq s_1 \leq ...\leq s_j \leq s,\; j \in \N^*, f \in {\cal C}_b(\R^k)\}
\label{eqW}
\end{equation}
where
\begin{equation}
f_{s_1,..,s_j}(\omega)=f((X_{s_1}(\omega),...,X_{s_j}(\omega))
\label{eqW2}
\end{equation}
\label{notcor2}
\end{notation}
From Proposition \ref{propDCC2}, we deduce the following Corollary:
\begin{corollary}
Assume that $\Theta$ satisfied hypothesis $H_{\theta}$ and has the Feller property.
 For all $Y \in {\cal W}^r_t$, $T^{\theta}_s(Y)$ belongs to ${\cal W}^r_s$.
\label{cor2}
\end{corollary}

\section{Time Consistent Convex Continuous Dynamic Procedures} 
\label{secTCC}
In this Section the Polish space is either 
${\cal C}^r$ or ${\cal D}^r$. We assume that  $\Theta$ satisfies hypothesis $H_{\theta}$ and has the Feller property.
The goal of this Section is to construct  time consistent convex procedures having a Feller  property, associated to a stable set of probability measures ${\cal Q}^{\Theta}(\Gamma)_{r,y}$ generated by $\Theta (\Gamma)$ (cf Corollary \ref{corsta}). 
%We start with a result that we prove for every separable metrizable space $Omega$. In particular it is satisfied for every Polish space and thus for ${\cal C}^r$ and  ${\cal D}^r$.
\begin{proposition}
Let $F=\{f_{\alpha},\; \alpha \in A\}$ be a set of l.s.c. functions uniformly bounded from below on $(\R^n)^k$.
\begin{itemize}
\item
There is a countable subset $\{f_{\alpha_n},\; n \in \N\}$ of $F$ such that 
\begin{equation}
\sup_{n \in \N}f_{\alpha_n}=\sup_{\alpha \in A}f_{\alpha},\; 
\label{eqqs}
\end{equation}
\item If the set $F$ is a lattice upward directed, there is an increasing sequence $g_n$ of elements of $F$ such that 
\begin{equation}
\sup_{\alpha \in A}f_{\alpha}=\lim_{n \rightarrow \infty}g_n,\; 
\label{eqqs2}
\end{equation}
\end{itemize}
\label{lemmaqs}
\end{proposition}
{\bf Proof}
\begin{enumerate}
\item Let $(K_p)_{p \in \N}$ be the closed ball of radius $p$ in $(\R^n)^k$ centered in $0$.  Let $(\phi_p)_{p \in \N^*}$ be an increasing sequence of continuous functions with compact support contained in $K_p$ such that the restriction of $\phi_p$ to $K_{p-1}$ is equal to $1$.
\item  Every function l.s.c. bounded from below is the limit of an increasing sequence of continuous bounded functions uniformly bounded from below. It follows then from {\it 1} that for all $\alpha$, there is a sequence $g^{\alpha}_p$ of continuous functions with compact support contained in $K_p$ such that $f_{\alpha}$ is the increasing limit of $g^{\alpha}_p$. Thus
\begin{equation}
\sup_{\alpha \in A} f_{\alpha}=\sup_{\alpha \in A,\; p \in \N^*} g^{\alpha}_p
\label{eqsci}
\end{equation}
 The set of continuous functions with compact support contained in $K_p$ is metrizable separable. It follows that the set $\{g^{\alpha}_p,\;\alpha \in A,\; p \in \N^*\}$ admits a countable dense subset.  Thus  there is a sequence $\alpha_n$ such that  
\begin{equation}
\sup_{\alpha \in A,\; p \in \N^*} g^{\alpha}_p=\sup_{n \in \N,\; p \in \N^*} g^{\alpha_n}_p
\label{eqsci2}
\end{equation}
By construction, for all $n,p$, $g^{\alpha_n}_p \leq f_{\alpha_n}$. Thus equation (\ref{eqqs}) follows from  (\ref{eqsci}) and (\ref{eqsci2}).
\item 
Let $f_{\alpha_n}$ in $F$ be  such that equation (\ref{eqqs})  is satisfied. Define the sequence $g_n$ by iteration: $g_0=f_{\alpha_0}$ and choose for $g_{n+1}$ a function  in $F$ such that $g_{n+1}\geq \sup(g_n,f_{\alpha_{n+1}})$. The sequence $g_n$ is by construction increasing. Let  $g$ its limit,  $\sup_{\alpha \in A}f_{\alpha}\geq g \geq \sup_{n \in \N} f_{\alpha_n}$. The result follows from (\ref{eqqs}). 
\end{enumerate} \hfill $\square $\\
We want now to construct general time consistent continuous procedures on ${\cal C}^r$ or ${\cal D}^r$.
For this we introduce the notion of Feller property for a penalty. Recall that a penalty has been defined  in Definition \ref{defpen}.
\begin{definition}
 Let ${\cal Q}={\cal Q}^{\Theta}(\Gamma)_{r,y}$ be a stable set of probability measures on $(\Omega^r,{\cal B}^r_T)$ generated by $\Theta (\Gamma)$ (cf Corollary \ref{corsta}). A penalty $(\alpha_{s,t})_{r \leq s \leq t\leq T}$ defined on ${\cal Q}$ is a  Feller penalty if 
\begin{equation}
\forall \theta \in \Theta(\Gamma),\; \exists \beta^{\theta}_{st} \in {\cal C}(\R^n),\; \alpha_{s,t}(Q^{\theta}_{r,y})= \beta^{\theta}_{st}(X_s)
\label{eqpenfe}
\end{equation}
\label{defpf}
\end{definition}
\begin{remark}
The above Feller property is required only for the probability measures $Q^{\theta}_{r,y}$ generating the stable set ${\cal Q}$ but not for all the probability measures $Q^{\gamma}_{r,y}$ ($\gamma \in \Theta(\Gamma)^r_T$) in ${\cal Q}$.  From the definition of a penalty (Definition \ref{defpen}), it follows that the function $\beta^{\theta}_{st}$ is bounded from above. We do not ask for boundedness.
\label{remFell}
\end{remark}

\begin{notation}
Let $0 \leq r \leq t$.  Recall that ${\cal W}^r_t$ (Notation \ref{notcor2}) is the algebra of continuous coordinates functions. 
%and ${\cal V}^r_t$ (Notation \ref{notcor})is the vector space defined from continuous functions with compact support. 
Denote also ${\cal L}_t$ the algebra ${\cal L}_t=\{f(X_t), \; f \text{continuous bounded on } \;\R^n \}$. In the following 
\begin{enumerate}
\item   ${\cal H}^r_t$ denotes either ${\cal W}^r_t$ or ${\cal L}_t$.
%\item In case $\Omega^r={\cal C}^r$,  ${\cal H}^r_t$ denotes either ${\cal W}^r_t$ or ${\cal V}^r_t$ or ${\cal L}_t$. 
\item   $\hat {\cal H}^r_t$ denotes either $\hat {\cal W}^r_t$ or $\hat {\cal L}_t$, with 
\begin{equation}
 \hat {\cal W}^r_t=\{f(X_{t_1},...X_{t_k}),\; r \leq t_1 <..<t_k\leq t,\;\; \text{f lsc bounded from below}\}
\label{eqwhat}
\end{equation}

\begin{equation}
\hat {\cal L}_t=\{f(X_t)\; \text{f lsc bounded from below}\}
\label{eqhhat}
\end{equation}
\end{enumerate}
\label{notH}
\end{notation}

The following proposition is an adaptation to this general martingale setting of the ideas used in the proof of Theorem  4.8 of \cite{JBN-PDE}.
\begin{proposition}
Let $r \leq s \leq t$. Assume that $\Theta$  satisfies hypothesis $H_{\theta}$ and the Feller property. Let $\Gamma$ be a  multivalued Borel mapping. Let  $\alpha_{st}$ be  a Feller penalty on  the  stable set of probability measures on $(\Omega^r,{\cal B}^r_T)$ generated by $\Theta (\Gamma)$. Let $Y \in \hat {\cal H}^r_t$, 
\begin{enumerate}
\item  For all $\theta \in \Theta(\Gamma)$,  
$T^{\theta}_s(Y)- \alpha_{s,t}(Q^{\theta}_{r,y})$ belongs to $ \hat {\cal H}^r_s$
\item
 For all $\gamma \in \Theta(\Gamma)^r_t$ and $s \in [r,t[$, there is $\eta_s \in  \Theta(\Gamma)^r_t$ such that 
\begin{equation}
T^{\gamma}_s(Y)- \alpha_{s,t}(Q^{\gamma}_{r,y})\leq T^{\eta_s }_s(Y)- \alpha_{s,t}(P^{\eta_s }_{r,y})=Z_s
\label{eqcontin}
\end{equation}
and such that $T^{\eta_s }_s(Y)- \alpha_{s,t}(P^{\eta_s }_{r,y})=Z_s$ belongs to $ \hat {\cal H}^r_s$.
\end{enumerate}
\label{propc}
\end{proposition}
{\bf Proof} 
Notice first that for all $\gamma \in \Theta(\Gamma)^r_t$ the map $T^{\gamma}_s$ is defined on all functions ${\cal B}^r_t$ measurable bounded from below and that if $Y,Z$ are ${\cal B}^r_t$ measurable bounded from below, $T^{\gamma}_s(Y+Z)=T^{\gamma}_s(Y)+T^{\gamma}_s(Z)$. This follows from the linearity of $T^{\gamma}_s$ on bounded functions and from the continuity from below.
\begin{enumerate}
\item Let $\theta \in \Theta(\Gamma)$. Let   $Y$ in  ${\cal H}^r_t$, it follows from Proposition  
%\ref{propCCCa} and 
\ref{propDCC2} that $T^{\theta}_s(Y)$ belongs to ${\cal H}^r_s$ . If  $Y$ belongs to $\hat {\cal H}^r_t$, $Y$ is the increasing limit of a sequence $Y_n \in {\cal H}^r_t$. From the continuity from below of $T^{\theta}_s$ (Corollary \ref{corext}), it follows that $T^{\theta}_s(Y)$ belongs to $ \hat {\cal H}^r_s$. It follows from the definition of a Feller penalty (Definitions \ref{defpf} and \ref{defpen}) that $-\alpha_{st}(Q^{\theta}_{r,y})$ belongs to $\hat {\cal H}^r_s$.  This gives the result.
\item
$r=s_0<s_1<...<s_n=t$ be a subdivision associated to $\gamma$. 
  For all $u \in [s_{n-1},s_n[$, $\gamma_{u}=\sum_{j \in I_{n-1}} 1_{A_{{n-1},j}}(\omega)\theta_{{n-1},j}(u,X_u(\omega))$. For all $j \in I_{n-1}$, and $s \in [s_{n-1},s_n[$, from 1., there is an element $Y_{s,j}$ of $ \hat {\cal H}^r_s$  such that $T^{\theta_{{n-1},j}}_{s}(Y)- \alpha_{s,t}(Q^{\theta_{{n-1},j}}_{r,y})=Y_{s,j}$. $Z_s=\sup_{j \in I_{n-1}}Y_{s,j}$ belongs to $ \hat {\cal H}^r_s$.  Let $B_{s,j}$ the ${\cal B}^r_s$ measurable  sets such that  $Z_s=\sum_{j \in I_{n-1}}1_{B_{s,j}}Y_{s,j}$. Let $\eta_s \in  \Theta(\Gamma)^r_t$ such that  $\eta_s(u)(\omega)=\sum_{j \in I_{n-1}}1_{B_{s,j}}(\omega)\theta_{{n-1},j}(u,X_u(\omega))$ for all $u \in [s,t]$. $\eta_s$ satisfies the required conditions.\\ 
\item For $s \in [s_i,s_{i+1}[$ with $i+1<n$, the result is  proved downward by induction.  $\gamma_{s}=\sum_{j \in I_{i}} 1_{A_{i,j}}(\omega)\theta_{i,j}(u,X_u(\omega))$ Assume that $Z_v$ and $\eta_v$ satifying equation(\ref{eqcontin}) and $Z_v \in \hat {\cal H}^r_s$ have been constructed  for all $v \geq s_{i+1}$. For all $j \in I_{i}$,  $T^{\theta_{i,j}}_{s}(Z_{s_{i+1}})- \alpha_{s,t}(Q^{\theta_{i,j}}_{r,y})=\phi_{s,j}$ belongs to $ \hat {\cal H}^r_s$. Let  $B_{s,j} \in {\cal B}^r_s$ such  that $Z_s=\sup_{j \in I_{i}}\phi_{s,j}=\sum_{j \in I_{i}}1_{B_{s,j}}\phi_{s,j}$. Let $\eta_s \in  \Theta(\Gamma)^r_t$ such that $\eta_s(u)=\eta_{s_{i+1}}(u)$ for all $u \geq s_{i+1}$ and $\eta_s(u)(\omega)=\sum_{j \in I_{i}}1_{B_{s,j}}\theta_{{i},j}(u,X_u(\omega))$ for all $u \in [s,s_{i+1}[$. $\eta_s$ and $Z_s$ satisfy the required conditions.
\end{enumerate}

\begin{theorem}
Let $r \leq s \leq t \leq T$. Let $y \in \R^n$. Assume that $\Theta$  satisfies hypothesis $H_{\theta}$ and the Feller property. Let $\Gamma$ be a  multivalued Borel mapping. 
Let  ${\cal Q}={\cal Q}^{\theta}(\Gamma)_{r,y}$ be the  stable set of probability measures generated by $\Theta(\Gamma)$. Let $(\alpha_{st})$ be a Feller penalty on ${\cal Q}$. Denote $T^{\gamma}_s=T^{Q^{\gamma}_{r,y}}_s$.
 For all $Y$ in $ \hat {\cal H}^r_t$, let 
\begin{equation}
\Pi^{r,y}_{s,t}(Y)=\sup_{Q^{\gamma}_{r,y} \in {\cal Q}^{\theta}(\Gamma)_{r,y}}  (T^{\gamma}_s(Y)- \alpha_{s,t}(Q^{\gamma}_{r,y}))
\label{eqtcd0}
\end{equation}
\begin{enumerate}
\item For all $Y \in  \hat {\cal H}^r_t$. 
$\Pi^{r,y}_{s,t}(Y)$ belongs to $ \hat {\cal H}^r_s$. Moreover for all $s \in [r,t]$, there is a sequence $Q_n=Q^{\gamma_n}_{r,y}$ in ${\cal Q}$, $\gamma_n$ in $\Theta(\Gamma)^r_t$  such that $T^{\gamma_n}_s(Y)-\alpha_{s,t}(Q^{\gamma_n}_{r,y})$ belongs to $\hat {\cal H}^r_s$, and such that $\Pi^{r,y}_{s,t}(Y)$ is the increasing limit:
\begin{equation}
\Pi^{r,y}_{s,t}(Y)=\lim  [T^{\gamma_n}_s(Y)- \alpha_{s,t}(Q^{\gamma_n}_{r,y})]
\label{eqtcd1}
\end{equation}

\item 
$\Pi^{r,y}_{s,t}$ is   a  convex   monotone map continuous from below on $ \hat {\cal H}^r_t$  with values in $ \hat {\cal H}^r_s$.\\

\item For all $r \leq s \leq t \leq u$, for all $Y \in \hat {\cal H}^r_u$,
\begin{equation}
\Pi^{r,y}_{s,u}(Y)=\Pi^{r,y}_{s,t}(\Pi^{r,y}_{t,u}(Y))
\label{tc510}
\end{equation}
\item In case ${\cal H}^r_t={\cal L}_t$, the sequence $\gamma_n$ of \it{1.} can be choosen independently of $r$.
\end{enumerate}
\label{ThmTCCP}
\end{theorem}
{\bf Proof} 
\begin{enumerate}
\item 
Let $Y \in \hat {\cal H}^r_t$.  For $Q^{\gamma}_{r,y}\in {\cal Q}$, let $Y^{\gamma}_s=T^{\gamma}_s(Y)- \alpha_{s,t}(Q^{\gamma}_{r,y})$. Let ${\cal Z}=\{Y^{\gamma}_s\; Q^{\gamma}_{r,y}\in {\cal Q},\; \gamma \in \Theta(\Gamma)^r_t\; |\;Y^{\gamma}_s \in \hat {\cal H}^r_s\}$. It follows from the stability by bifurcation of  ${\cal Q}$, from the expression of $T^{\gamma}$ (equation (\ref{Tgamma}))  and  from the locality of the penalty  that ${\cal Z}$ is a lattice upward directed.
\begin{itemize}
\item If ${\cal H}^r_t={\cal W}^r_t$, there is a lsc function $f$ bounded from below such that $Y=f(X_{s_1},X_{s_2},..X_{s_l})$. Let $j$ such that  $s_j < s \leq s_{j+1}$. From Proposition \ref {propDCC2} for all $\theta \in \Theta(\Gamma)$, $Y^{\theta}_s=g(X_{s_1},..., X_{s_j},X_s)$  for some  function $g$ continuous  bounded from below on $(\R^n)^{j+1}$. It follows from the proof of Proposition \ref{propc} that $\{g\;|g(X_{s_1},..., X_{s_j},X_s) \in {\cal Z}\}$ is a lattice of lsc  functions bounded from below  on $(\R^n)^{j+1}$. The result follows from Proposition \ref{lemmaqs}.

\item In  case  ${\cal H}^r_t={\cal L}_t$, $Y=f(X_t)$. For all $\gamma \in \Theta(\Gamma)^r_T$, $Y^{\gamma}_s=g(X_s)$  for some  function $g$ Borelian  bounded from below on $(\R^n)^{j+1}$ independent of $r$. A  similar proof as the above one proves {\it 1.} and {\it 4.}
\end{itemize}
\item From the linearity of $T^{\gamma}_s$ on bounded functions and the continuity form below, it follows that for all non negative real numbers $c$ and $d$, ,  for all $Y,Z$ in ${\cal H}^r_t$, $T^{\gamma}_s (cY+dZ)= c T^{\gamma}_s (Y)+dT^{\gamma}_s (Z)$. The convexity of $\Pi^{r,y}_{s,t}$ follows. Let $Y_n$ be increasing to $Y$. Every $T^{\gamma}_s$ is continuous from below  thus $\Pi^{r,y}_{s,t}( Y)= \sup_{\gamma}T^{\gamma}_s(\lim Y_n)=\sup_{\gamma,n}T^{\gamma}_s(Y_n)=\sup_n\Pi^{r,y}_{s,t}(Y_n)$.
\item Let $Y \in \hat {\cal H}^r_u$. From 1., $\Pi^{r,y}_{t,u}(Y)$  belongs to $\hat {\cal H}^r_t$ and we have the following  increasing limits:
\begin{equation}
\Pi^{r,y}_{s,t}(\Pi^{r,y}_{t,u}(Y))=\lim_{n \rightarrow \infty}  [T^{\gamma_n}_s(\Pi^{r,y}_{t,u}(Y))- \alpha_{s,t}(Q^{\gamma_n}_{r,y})]
\label{eqtcd2}
\end{equation}
 \begin{equation}
\Pi^{r,y}_{t,u}(Y)=\lim_{k \rightarrow \infty}  [T^{\delta_k}_t((Y))- \alpha_{t,u}(Q^{\delta_k}_{r,y})]
\label{eqtcd3}
\end{equation}
Fom the continuity from below of $T^{\gamma_n}_s$, it follows that 
\begin{equation}
\Pi^{r,y}_{s,t}(\Pi^{r,y}_{t,u}(Y))=\sup_{n,k}  (T^{\gamma_n}_s[T^{\delta_k}_t((Y))- \alpha_{t,u}(Q^{\delta_k}_{r,y})]- \alpha_{s,t}(Q^{\gamma_n}_{r,y}))
\label{eqtcd4}
\end{equation}
From the stability of ${\cal Q}$,  for all $n$ and $k$, let $\lambda_{n,k}$ be defined as in Theorem \ref{thmsta} i), such that $T^{\gamma_n}_s(\phi)=T^{\lambda_{n,k}}_s(\phi)$ for all $\phi$ ${\cal B}^r_t$ measurable bounded and $T^{\delta_k}_t=T^{\lambda_{n,k}}_t$. From the chain rule, $T^{\gamma_n}_s \circ [T^{\delta_k}_t]= T^{\lambda_{n,k}}_s$. From the definition of $\lambda_{n,k}$, the local property of the penalty and the cocycle condition, it follows that $T^{\gamma_n}_s[- \alpha_{t,u}(Q^{\delta_k}_{r,y})]- \alpha_{s,t}(Q^{\gamma_n}_{r,y})=- \alpha_{s,u}(Q^{\lambda_{n,k}}_{r,y})$. 
Thus from (\ref{eqtcd4}), it follows that 
\begin{equation}
\Pi^{r,y}_{s,t}(\Pi^{r,y}_{t,u}(Y))=\sup_{n,k}  [ T^{\lambda_{n,k}}_s(Y))- \alpha_{s,u}(Q^{\lambda_{n,k}}_{r,y})] \leq \Pi^{r,y}_{s,u}(Y)
\label{eqtcd5}
\end{equation}
Conversely, $\Pi^{r,y}_{s,u}(Y)$ is the increasing limit of $T^{\nu_j}_s(Y)- \alpha_{s,u}(Q^{\nu_j}_{r,y})$. Making use of the chain rule for $T^{\nu_j}$ and of the cocycle condition, we get the inequality  $\Pi^{r,y}_{s,u}(Y) \leq \Pi^{r,y}_{s,t}(\Pi^{r,y}_{t,u}(Y))$.
\item If $Y \in {\cal L}_t$, $Y=f(X_t)$. 
\end{enumerate}
\begin{definition}
Let $r \leq s \leq t \leq T$. Let $\Theta$ satisfying hypothesis $H_{\theta}$.
let  ${\cal Q}={\cal Q}^{\theta}(\Gamma)_{r,y}$ be the  stable set of probability measures generated by $\Theta(\Gamma)$. Let $(\alpha_{st})$ be a Feller penalty on ${\cal Q}$.
 For all $Y$ in $ \hat {\cal H}^r_t$, let 
\begin{equation}
\Pi^{r,y}_{s,t}(Y)=\sup_{Q^{\gamma}_{r,y} \in {\cal Q}^{\theta}(\Gamma)_{r,y}}  (T^{Q^{\gamma}_{r,y}}_s(Y)- \alpha_{s,t}(Q^{\gamma}_{r,y})
\label{eqtcd6}
\end{equation}
One says that $\Pi^{r,y}_{s,t}$ is a time consistent dynamic procedure on $ \hat {\cal H}^r_t$ if all properties  of Theorem \ref{ThmTCCP} are satisfied.
\label{defTCPP}
\end{definition}
\begin{proposition} Let ${\cal P}={\cal Q}^{\theta}(\Gamma)_{r,y}$. Denote $\overline {\cal W}^r_t$ the closure of the lattice vector space ${\cal W}^r_t$ for the norm $\sup_{Q \in {\cal P}}E_Q(|X|)$. Assume that for all $s$, and $X \in {\cal W}^r_t$ , $\Pi^{r,y}_{s,t}(X)$ belongs to $\overline {\cal W}^r_s$. The procedure $\Pi^{r,y}_{s,t}$can be uniquely extended to $\overline {\cal W}^r_t$. The time consistency (equation \ref{tc510}) extends to  $\overline {\cal W}^r_t$.\\
In case where ${\cal P}$ is a weakly relatively compact set of probability measures on ${\cal C}([r,T],\R^n)$, $\overline {\cal W}^r_t$ is equal to $L^1(c)$.
\end{proposition}
{\bf Proof}
Let $X,Y$ in ${\cal W}^r_t$. From equation (\ref{eqtcd6}), 
$\Pi^{r,y}_{s,t}(Y)$ is the increasing limit of $T^{\gamma_n}_s(Y)- \alpha_{s,t}(Q^{\gamma_n}_{r,y})$. It follows that $\Pi^{r,y}_{st}(Y) \leq \Pi^{r,y}_{st}(X)+ sup_{n \in \N}T^{\gamma_n}_s(Y-X)$. Changing the roles of $X$ and $Y$, it follows that 
\begin{equation}
|\Pi^{r,y}_{st}(Y)-\Pi^{r,y}_{st}(X)| \leq \sup_{n,k}[T^{\gamma_n}_s(|Y-X|), T^{\delta_k}_s(|Y-X|)]
\label{eqextens}
\end{equation}
 As in the proof of equation (\ref{eqtcd1}), it follows from Proposition \ref{lemmaqs} that \\$\sup_{Q^{\gamma} \in {\cal P}}T^{\gamma}_s(|X-Y|)$ is the increasing limit of a sequence $T^{\alpha_n}_s (|X-Y|)$. It follows from equation (\ref {eqextens})  and the stability of ${\cal P}$ that 
\begin{equation}
\sup_{P \in {\cal P}} E_P(|\Pi^{r,y}_{st}(Y)-\Pi^{r,y}_{st}(X)|) \leq 
\sup_{P \in {\cal P}} E_P(|X-Y|)
\label{eqL1C}
\end{equation} 
 It follows from equation (\ref{eqL1C}), that $\Pi^{r,y}_{st}$ can be uniquely extended to $\overline {\cal W}^r_t$, with values in  $\overline {\cal W}^r_s$, equation (\ref{tc510}) is then satisfied on $\overline {\cal W}^r_t$.\\
The last result follows from Lemma \ref{lemmacont}.

\section{ Construction of Feller penalties}
\label{secTCCC}
\subsection{Penalties for continuous diffusions}
\label{secTCCC0}
In this Section we restrict to the case of continuous paths. We consider the martingale problem for diffusions.   In \cite{JBN-PDE}  we have constructed Feller penalties associated to probability measures solution to the martingale problem for continuous diffusions. In \cite{JBN-PDE} the probability measures in the stable set were all solutions to a martingale problem with the same function $a(t,x)$. In the present paper the probability measures are no more equivalent, they are even not dominated. The construction of the penalties that we describe below  is the generalization of the construction that we have done in \cite{JBN-PDE}. The ideas are the same.
In this paper we chose the parameters $a$ and $b$.

\begin{definition} 
$\Theta$ is the set of continuous bounded functions $\theta(t,x)=$\\ $(a(t,x),b(t,x)))$ with values in $M_n(\R)) \times \R^n$ such that for all $(t,x)$ the matrix $a(t,x)$ is invertible.
Let $\Gamma$ be  a multivalued Borel mapping from $\R_+ \times \R^n$ to $M_n(\R) \times \R^n$, $\Theta(\Gamma)=\{\theta \in {\Theta},\; \Gamma\;\text{valued}\}$. 
\label{deftetacont}
\end{definition}
For general $\gamma=(\eta,\mu) \in {\Theta}^r_T$  we denote $P^{\eta,\mu}_{r,y}$ the unique solution to the martingale problem starting from $y$ at time $r$ constructed previously. In the case  $\gamma \in \Theta$, $\gamma(u,\omega)=\theta(u,X_u(\omega))$ for some continuous $\theta(u,x)$. When we want to emphasize the fact that we refer to  continuous parameters $\eta$ and $\mu$, the probability measure $P^{\eta,\mu}_{r,y}$ will be denoted $Q^{a,b}_{r,y}$.

As in \cite{JBN-PDE} we introduce hypothesis $H_g$.
\begin{definition}
Hypothesis $H_g$\\
\begin{enumerate}
\item 
$g:\R_+ \times \R^n \times (M_n(\R) \times \R^n) \rightarrow \R\cup {\infty}$ is a ``Caratheodory function on $\Gamma$''\\
More precisely, $g$ is Borelian and for all $u$, the restrition of $g_u$ to $\{(x,y), y \in \Gamma(u,x)\}$ is continuous  ($g_u(x,y)=g(u,x,y)$).
\item $g$ has polynomial growth on $\Gamma$ 
\begin{equation}
\forall (t,x,y), \; y  \in \Gamma(t,x),\;\;|g(t,x,y)| \leq C(1+||x||^m)
\label{eqpolg}
\end{equation}
\item $g$ is bounded from above on $\Gamma$ which means that $g$ is bounded from above on $\{(t,x,y), \; y  \in \Gamma(t,x)\}$
\end{enumerate}
\label{defHgc}
\end{definition}
This last condition on $\Gamma$ which was not present in \cite{JBN-PDE} could be suppressed with a little more work.\\
 Recall that from Corollary \ref{corsta}, the set ${\cal Q}^{\Theta}(\Gamma)_{r,y}=\{Q^{\gamma}_{r,y},\; \gamma \in {\Theta}(\Gamma)^r_T\}$ is stable. 
\begin{proposition} Let $\Gamma$ be a multivalued Borel mapping. Let 
 $\Theta(\Gamma)$ be as in Definition \ref{deftetacont}.   Assume that $g$ satisfies Hypothesis $H_g$. Let $t>0$.
\begin{enumerate}
\item  For all $\theta=(a,b) \in \Theta(\Gamma)$, there is a real valued continuous  map $(s,x) \in [0,t]\times \R^n \rightarrow L^{a,b}_t(g)(s,x)$  such that 
\begin{equation}
Q^{a,b}_{s,y}(\int_s^t g(u, X_u,\theta (u,X_u)) du )=L^{a,b}_t(g)(s,y) \;\;\;\forall s \in [0,t]\;\;and\; y \in \R^n
\label{eqFell47-3}
\end{equation}
\item For all $0 \leq r \leq s \leq  t$ and all $y \in \R^n$, for $Q=Q^{a,b}_{r,y}$, 
\begin{equation}
T^{Q}_s[-\int_s^t g(u, X_u,\theta(u,X_u)) du ]=-L^{a,b}_t(g)(s,X_s(\omega))
\label{eqFell47b-3}
\end{equation}
\end{enumerate}
Notice that for all $y \in \R^n$ $L^{a,b}_t(g)(t,x)=0$
\label{propcont3}
\end{proposition}  
{\bf Proof} {\it 1.} is an application of the first statement of Proposition 4.2 of \cite{JBN-PDE}\\
{\it 2.} $-g$ is bounded from below, so making use of the extension of $T^{Q}_s$ defined on Corollary \ref{corext} it follows  that 
\begin{eqnarray}
T^{Q}_s[-\int_s^t g(u, X_u,\theta(u,X_u)) du ]&=&Q^{a,b}_{s,X_s(\omega)}[-\int_s^t g(u, X_u,\theta(u,X_u)) du ]\nonumber\\
&=&- L^{a,b}_t(g)(s,X_s(\omega))
\label{eqfepe}
\end{eqnarray}
We  define now the penalty of every  probability measure $Q^{\gamma}_{r,y}$ in ${\cal Q}^{\Theta}(\Gamma)_{r,y}$. 
\begin{definition} For all $Q^{\gamma}_{r,y}$, for all $r \leq s\leq t \leq T$, we define the penalty 
\begin{equation}
\alpha_{s,t}(Q^{\gamma}_{r,y})=-T^{\gamma}_s[-\int\limits_{s}^{t} g(u,X_u(\omega),\gamma(u,\omega)) d u]
\label{eqP1_0}
\end{equation}
\label{defpen0}
\end{definition}

\begin{corollary}
The above definition provides a Feller penalty according to Definitions \ref{defpen} and  \ref{defpf}.
\label{corpen}
\end{corollary}
%\begin{corollary}
{\bf Proof} It follows from hypothesis $H_g$ that for all $r \leq s \leq t$, the function $-\int\limits_{s}^{t} g(u,X_u(\omega),\gamma(u,\omega)) d u]$ is ${\cal B}^r_t$ measurable bounded from below, thus from Corollary \ref{corext}, 
$\alpha_{s,t}(Q^{\gamma}_{r,y})$ is ${\cal B}^r_s$ measurable bounded from above.  The cocycle condition for $\alpha_{s,t}(Q^{\gamma}_{r,y})$ follows from the chain rule for $T^{\gamma}$ (  Corollary \ref{corext})and the additivity of $T^{\gamma}_s$ on functions bounded from below. 
The local condition follows from the Definition of the penalty. The Feller property of the penalty follows from  Proposition \ref{propcont3}.

\subsection{Penalties for  diffusions with Levy generators}
\label{secpenlevy}
In this Section we consider probability measures on the Polish space  ${\cal D}^r$ of c\`adl\`ag paths endowed with the Skorhokod topology. The notations are those of Section \ref{secDiLev}.
\begin{definition} 
\begin{enumerate}
\item
Let $C>0$. ${\cal M}_C$ denotes the set of $\sigma$-finite measures $\mu$ on $\R^n -\{0\}$ such that 
\begin{equation}
\int [||y||^2 1_{||y|| \leq 1}+||y||1_{||y|| >1}]\mu(dy) \leq C
\label{eqMc}
\end{equation}
\item 
$\Theta$ is the set of $(a,b,M)$ such that
\begin{itemize}
\item
 $a,b$ are continuous bounded functions on $\R_+ \times \R^n$, with values in $M^n(\R)$, $\R^n$, and for all $(t,x)$ $a(t,x)$ is invertible.
\item  $M$ satisfies hypothesis $M_C$ (Definition \ref{defMc}), i.e. for all $(t,x)$, $M(t,x) \in {\cal M}_c$. Furthermore for all Borelian subset 
$\Delta$ of $\R^n-\{0\}$,
\begin{equation} 
\int _{\Delta}\frac{||y||^2}{1+||y||^2} M(s,x,dy)\;
\label{eqbound1}
\end{equation}
is a continuous bounded function of $(s,x)$.
\end{itemize}
\item
${\cal M}_c$ is endowed with the weak topology. Let $\Gamma$ be  a closed convex multivalued Borel mapping from $\R_+ \times \R^n$ to $M_n(\R) \times \R^n \times {\cal M}_c$. Let $\Theta(\Gamma)=\{\theta \in {\Theta},\; \Gamma\;\text{valued}\}$. 
\end{enumerate}
\label{deftetacont2}
\end{definition}
We  want now to define  the penalty of every  probability measure $Q^{\gamma}_{r,y}$ in ${\cal Q}^{\Theta}(\Gamma)_{r,y}$. 
Recall that $\Theta^r_T$ is described by Definition \ref{Slambda}. For all $\gamma \in \Theta(\Gamma)^r_T$, write $\gamma(u,\omega)=(\eta,\lambda,\mu)(u,\omega)$,  where $\eta$ takes values in $M_n(\R)$, $\lambda$ in $\R^n$ and $\mu$ takes values in ${\cal M}_C$.
\begin{definition} 
Let $g$ be a real valued function on $\R^+ \times M_n(\R)\times \R^n \times {\cal M}_C$ such that for all $\theta=(a,b,M) \in \Theta(\Gamma)$, $ g^{\theta}(u,x)=g(u,x,a(u,x),b(u,x),M(u,x)$ is continuous bounded.
We define the penalty 
\begin{equation}
\alpha_{s,t}(Q^{\gamma}_{r,y})=T^{\gamma}_s[\int\limits_{s}^{t} g(u,X_u(\omega),\eta(u,\omega),\lambda(u,\omega),\mu(u,\omega))d u]
\label{eqP1_}
\end{equation}
\label{defpenlev}
\end{definition}
%Notice that from hypothesis $(M)$, for all Borelian subset $\Delta$ of $\R^n-\{0\}$, $1_{\Delta}\frac{||y||^2}{1+||y||^2}$ is a possible choice for the function $\phi$. 
\begin{proposition}
The above definition provides a penalty according to Definition \ref{defpen}.
 The penalty is a Feller penalty (Definition \ref{defpf}). More precisely, for all $\theta=(a,b,M) \in \Theta(\Gamma)$, there is a continuous function $L_{t}^{\theta}(g)$ on $[0,t] \times \R^n$ such that $\forall s \in [0,t]\;\;and\; y \in \R^n$,
\begin{equation}
Q^{a,b,M}_{s,y}(\int_s^t g(u,X_u(\omega),a(u,X_u(\omega)),b(u,X_u(\omega)), M(u,X_u(\omega))du= L_{t}^{\theta}(g)(s,y) \;\;\;
\label{eqFell47-4}
\end{equation}
\item For all $0 \leq r \leq s \leq  t$ and all $y \in \R^n$, for $Q=Q^{a,b,M}_{r,y}$, 
\begin{equation}
T^{Q}_s[\int_s^t g(u,X_u(\omega),a(u,X_u(\omega)),b(u,X_u(\omega)), M(u,X_u(\omega))du]=L_{t}^{\theta}(g)(s,X_s(\omega))
\label{eqFell47b-4}
\end{equation}
\label{proppen2}
\end{proposition}
{\bf Proof} 
\begin{enumerate}
\item It follows easily from the hypothesis, and the description of $\Theta(\Gamma)^r_T$ (cf Definition \ref{Slambda})  that for all $\gamma=(\eta,\lambda,\mu) \in \Theta(\Gamma)^r_T$, the function $\int\limits_{s}^{t} g(u,X_u(\omega),\eta(u,\omega),\lambda(u,\omega),\mu(u,\omega))$ is ${\cal B}^r_t$- measurable bounded. From Proposition \ref{propcanreg} it follows that $\alpha_{s,t}(Q^{\gamma}_{r,y})$ is  ${\cal B}^r_s$ measurable bounded.
\item The cocycle condition follows easily from the definition of the penalty (\ref{eqP1_}),  the linearity for $T^{\gamma}$ and the chain rule for $T^{\gamma}$.
\item 
locality  Let  $\gamma,\eta \in {\Theta^r_T}$, for all $A \in {\cal B}^r_s$, if $1_A \gamma(u,\omega)=1_A \eta(u,\omega)$, for all $u \in]s,t]$. from Proposition \ref{propcanreg} equation (\ref{Tgamma}), $T^{\gamma}_s$ and $T^{\eta}_s$ coincide on bounded variables ${\cal B}^r_t$ measurable. The locality follows then from the definition of the penalty. 
 \item Let $\theta \in \Theta(\Gamma)$, $\theta=(a,b,M)$, $a,b$ are continuous bounded functions of  $(u,x)$, and $M$ satisfies hypothesis $M_C$.  Let $\epsilon>0$. Let $K$   such that equation (\ref{eqLM0}) is  satisfied for all $0 \leq r \leq t$ and $||y|| \leq D$. The function $g^{\theta}$ is  bounded and uniformly  continuous   on $[0,T] \times \{||x|| \leq K\}$.  Making use of  equation (\ref{eqLM}) and of arguments similar to those of Proposition \ref{lemmFell},  it follows   that for all $s' \leq t$, there  are $s' \leq t_1<t_2<...<t_k \leq t$ and a continuous bounded function $f_{s',t}$ on $(\R^n)^k$ such that for all $||y|| \leq D$, for all $0 \leq r \leq s'$,
\begin{equation}
Q^{a,b,M}_{r,y}(|\int\limits_{s'}^{t} g^{\theta}(u,X_u(\omega))du-f_{s',t}(X_{t_1},X_{t_2},..X_{t_k})|<\epsilon
\label{equnif00}
\end{equation}
From Proposition  \ref{lemmFell}, there is a continuous function $\hat f_{s',t}$ on $[0,s'] \times \R^n$ such that $Q^{a,b,M}_{r,y}(f_{s',t}(X_{t_1},X_{t_2},..X_{t_k}))=\hat f_{s',t}(r,y)$. It follows from (\ref{equnif00}) that  $Q^{a,b,M}_{r,y}(\int\limits_{s}^{t} g^{\theta}(u,X_u(\omega))du$ is the uniform limit on $[0,s'] \times \{||y|| \leq D\}$ of a sequence of continuous functions. It is thus continuous. The function $g$ being bounded, the continuity on $[0,t] \times \R^n$ of the function $L^{\theta}_t$ defined by equation (\ref{eqFell47-4}) follows then easily.  Equation (\ref{eqFell47b-4}) follows  from Proposition \ref{proCRC}.
 \end{enumerate}
\section{Feller property of the time consistent dynamic procedure}
In all the following we restrict to the case of continuous diffusions or diffusions with Levy generator. The notations and hypothesis are those of Section \ref{secTCCC}. 
%$\Theta$ is either the set defined in  Definition \ref{deftetacont} in the case of continuous diffusions or the set defined in Definition \ref{deftetacont2} in the case of diffusions with Levy generator, assuming furthermore that for all $\theta=(a,b,M) \in \Theta$, $M$ satisfies hypothesis $M_C$ for some $C>0$.\\
 %Let $\Gamma$ be a  multivalued Borel mapping. Let  ${\cal Q}={\cal Q}^{\theta}(\Gamma)_{r,y}$ be the  stable set of probability measures generated by $\Theta(\Gamma)$. Let $(\alpha_{st})$ be the Feller penalty on ${\cal Q}$ constructed from $g$ satisfying hypothesis $H_g$ in Section \ref{secTCCC} for continuous diffusions and constructed as in Definition \ref{defpenlev} in Section \ref{secpenlevy} for diffusions with Levy generator. 
  We prove now that  with the specific choice for the penalty which  was made in Subsections \ref{secTCCC0} and \ref{secpenlevy}, the procedures $\Pi^{r,y}_{s,t}$ that we have constructed for all given $r$ and $y$ are strongly connected and have a Feller property.
We introduce now a larger set of parameters: $\tilde \Theta(\Gamma)^r_T$. The definition of $\tilde \Theta(\Gamma)^r_T$ is similar to  the definition of $\Theta(\Gamma)^r_T$ (cf Definition \ref{Slambda}), but  now  the set $I_0$ is  a finite set and not necessarily a singleton.  Every element $\tilde \gamma$  of $\tilde \Theta(\Gamma)^r_T$, satisfies  equation (\ref{eqstable}). For $\tilde \gamma \in \tilde \Theta(\Gamma)^r_T$, the process $(Z^{\tilde \gamma}_{r,u})_{r \leq u \leq T}$ is defined by equation (\ref{eqstable2}).
\begin{remark}
Let $0 \leq s <t$. Let $x \in \R^n$. Let $\tilde \gamma \in \tilde \Theta(\Gamma)^s_t$. The sets $A_{0,j}, j \in I_0$ belong to the $\sigma$-algebra generated by $X_s$ and form a partition of $\Omega^r_t$. It follows that there is $j_x \in I_0$ such that $\{X_s(\omega)=x\} \subset A_{0,j_x}$. Therefore a probability measure $Q$ on $(\Omega^s,{\cal B}^s_t)$ is solution to the martingale problem $(Z^{\tilde \gamma}_{s,u})_{s \leq u \leq t}$ starting from $x$ at time  $s$ (Definition \ref{defFellprop}) if and only if $Q$ is solution to the martingale problem $(Z^{\gamma_x}_{s,u})_{s \leq u \leq t}$ starting from $x$ at time  $s$, where $(\gamma_x)(u,\omega)=\theta_{0,j_x}(u,X_u(\omega))$ for all $s \leq u < s_1$, and $(\gamma_x) (u,\omega)=\tilde\gamma(u,\omega)$ for all $s_1 \leq u$. Thus ${\gamma_x}$ belongs to $\Theta(\Gamma)^s_t$.
\label{remtheta0}
\end{remark}

The following proposition  is a more precise version of Proposition \ref{propc}.
\begin{proposition}
Let $h$ be a continuous  function on $\R^n$ bounded from below.
 For all $\gamma \in \Theta(\Gamma)^r_t$ there is a continuous  function $h^{\gamma}$ on $[0,t] \times \R^n$ bounded from below and for all   $r \leq  s \leq t$, there is $\eta_s \in  \tilde \Theta(\Gamma)^s_t$ such that 
\begin{enumerate}
\item
\begin{equation}
\forall x \in \R^n,\;T^{P^{\eta_s}_{s,x}}_s(h(X_t))-\alpha_{s,t}(P^{\eta_s}_{s,x})=h^{\gamma}(s,x)
\label{eqetas}
\end{equation}
There is a partition $B_j$ of $\R^n$ in Borelian sets, and $\eta^j_s \in \Theta(\Gamma)^s_t$ such that
\begin{equation}
\forall x \in B_j,\;T^{P^{\eta^j_s}_{s,x}}_s(h(X_t))-\alpha_{s,t}(P^{\eta^j_s}_{s,x})=h^{\gamma}(s,x)
\label{eqetasj}
\end{equation}
\begin{equation}
\forall x \in B_j,\;h^{\eta^j_s}(u,x)=h^{\gamma}(u,x), \; \forall s \leq u \leq t
\label{eqetasj2}
\end{equation}
\item
For all $\eta \in \theta(\Gamma)^r_t$ whose restriction to $[s,t]$ is equal to $\eta_s$,
\begin{eqnarray}
\forall y  \in \R^n,\;\;T^{P^{  \eta}_{r,y}}_s (h(X_t)-\alpha_{s,t}(P^{\eta}_{r,y})=h^{\gamma}(s,X_s)\nonumber\\
T^{Q^{\gamma}_{r,y}}_s(h(X_t))-\alpha_{s,t}(Q^{\gamma}_{r,y}) \leq T^{P^{\eta}_{r,y}}_s(h(X_t))-\alpha_{s,t}(P^{ \eta}_{r,y})\nonumber\\
h^{\gamma}(u,x)=h^{\eta}(u,x) \;\;\forall s \leq u \leq t
\label{eqnfell}
\end{eqnarray}
\end{enumerate}
\label{propsup}
\end{proposition}
{\bf Proof}
\begin{itemize}
\item Step 1:
Let $\theta \in \Theta(\Gamma)$, let $h^{\theta}$ be the function defined  on $[0,t] \times \R^n$ by 
\begin{equation}
\;T^{P^{\theta}_{s,x}}_s(h(X_t))-\alpha_{s,t}(P^{\theta}_{s,x})=h^{\theta}(s,x)
\label{eqhteta}
\end{equation}
 The continuity of the map $h^{\theta}$ on $[0,t] \times \R^n$ follows from Theorem 7.1 of \cite{SV1} (see  Proposition 2.7 of \cite{JBN-PDE} for details), and from Proposition  \ref{propcont3} in case of diffusions.   
In case of diffusions with Levy generator, the continuity of $h^{\theta}$  follows from Propositions \ref{lemmFell} and  \ref{proppen2}. Equation 
\begin{equation}
\forall 0 \leq r \leq s \leq t, \forall y  \in \R^n,\;\;T^{P^{\theta}_{r,y}}_s (h(X_t))-\alpha_{s,t}(P^{\theta}_{r,y})=h^{\theta}(s,X_s)
\end{equation}
follows then from Lemma \ref{lemmaCC2}, Propositions  \ref{propcont3} and \ref{proppen2}. \\
\item Step 2:
The proof of the existence of $h^{\gamma}$ and $\eta_s$  satisfying   (\ref{eqnfell}) follows  the proof of
Theorem 4.8 in \cite{JBN-PDE}. We sketch the proof.
Let $r=s_0<s_1<...<s_n=t$ be a subdivision associated to $\gamma$. For  $u \in [s_i,s_{i+1}[$,  $\gamma(u,\omega)=\sum_{j \in I_i} 1_{A_{i,j}}(\omega)\theta_{i,j}(u,X_u(\omega))$, $\theta_{i,j} \in {\Theta}(\Gamma)$. $h^{\gamma}$ and $\eta_s$ are defined recursively. Assume that $h^{\gamma}$ is defined and continuous on $[s_{i+1},t] \times \R^n$, and that $\eta_u$ is defined for $u \in [s_{i+1},t]$.  Let $h_{i,j}$ be the continuous function on $[s_i,s_{i+1}]$ such that $\forall s \in [s_i,s_{i+1}]$,
\begin{equation}
\forall x \in \R^n,\;T^{P^{\theta_{i, j}}_{s,x}}_s(h^{\gamma}(s_{i+1},X_{s_{i+1}}))-\alpha_{s,s_{i+1}}(P^{\theta_{i,j}}_{s,x})=h_{i,j}(s,x)
\label{eqhteta2}
\end{equation}
For $s_i \leq s < s_{i+1}$, let $h^{\gamma}(s,x) =\sup_{j \in I_i} h_{i,j}(s,x)$. Given $s \in [s_i, s_{i+1}[$, there is a Borelian partition $B_{j},j \in I_i$  of $\R^n$ such that $h^{\gamma}(s,x) =\sum _{j \in I_i} 1_{B_{j}}(x) h_{i,j}(s,x)$. Let 
$\eta_s(u, \omega)=\eta^j_s(u, \omega)=\eta_{s_{i+1}}(u, \omega)$ for $u \geq  s_{i+1}$.  Let  $\eta^j_s(u,\omega)=\theta_{i,j}(u,X_u(\omega))$ and $\eta_s(u,\omega)= \sum_j1_{B_{j}}(X_s)\theta_{i,j}(u,X_u(\omega))$ for $s \leq u <s_{i+1}$. 
 We end the proof of {\it 1.} making use of the expression of $T^{Q^{\gamma}_{r,y}}_s$, cf  Theorem \ref{propcanreg} and of Remark \ref{remtheta0}. Notice that  $\eta_r$ belongs to $\Theta(\Gamma)^r_t$, $\eta^j_s \in  \Theta(\Gamma)^s_t$ and  $\eta_s \in  \tilde \Theta(\Gamma)^s_t$
\item Step 3:
It follows from Step 1, from the expressions of $T^{Q^{\gamma}_{r,y}}_s$,  $T^{P^{\eta}_{r,y}}_s$, and the properties of the penalty that the construction given above leads to $h^{\gamma}$ and $\eta_s$ such that  equation (\ref{eqnfell}) is satisfied for all  $\eta \in \theta(\Gamma)^r_t$ whose restriction to $[s,t]$ is equal to $\eta_s$. 
\end{itemize}\hfill $\square $ 
\begin{theorem}
\begin{enumerate}
\item
Let $0 \leq r \leq T$. Let $\hat {\cal W}^r_T$ be the cone  of  coordinate functions $f(X_{t_1},..X_{t_k})$, $f$ lower semi continuous bounded from below (Notation \ref{notH}). Let $(\Pi^{r,y}_{s,t})_{r \leq s \leq t \leq T}$ be defined on $\hat {\cal W}^r_T$ by  
\begin{equation}
\Pi^{r,y}_{s,t}(Y)=\sup_{Q^{\gamma}_{r,y} \in {\cal Q}^{\theta}(\Gamma)_{r,y}}  (T^{Q^{\gamma}_{r,y}}_s(Y)- \alpha_{s,t}(Q^{\gamma}_{r,y}))
\label{eqtcd7}
\end{equation}
where $\alpha_{s,t}$ is the penalty constructed in Section \ref{secTCCC}.
Then  $(\Pi^{r,y}_{s,t})_{r \leq s \leq t \leq T}$ is a  time  consistent convex  dynamic procedure on $\hat {\cal W}^r_t$. 
\item
It has the following Feller property: For all  $h$  lower semi continuous  function on $\R^n$ bounded from below, 
 there  is a lower semi-continuous  function $\tilde h$ bounded from below on   $[0,t] \times \R^n$  such that 
\begin{equation}
 \forall x \in \R^n,\;\Pi^{s,x}_{s,t}(h(X_t))=\tilde h (s,x)
\label{eqmeas20}
\end{equation}
\begin{equation}
\forall 0 \leq r < s \leq t, \;\forall  y \in \R^n,\;\;\Pi^{r,y}_{s,t}(h(X_t))=\tilde h (s,X_{s})
\label{eqmeas0}
\end{equation}
Furthermore for all given $r$ and $y$,  for all $r<s \leq t$, there is a sequence $\eta_n \in \Theta(\Gamma)^r_t$ such that $\Pi^{r,y}_{s,t}(h(X_t))$ is the increasing limit of $(T^{P^{\eta_n}_{r,y}}_s(h(X_t))-\alpha_{s,t}(P^{\eta_n}_{r,y}))$
\end{enumerate}
\label{thmfellerpro}
\end{theorem}
{\bf Proof} 
It follows from Corollary \ref {corpen} and Proposition \ref{proppen2} that $\alpha$ is a Feller penalty. Thus from Theorem \ref{ThmTCCP}, and Definition \ref{defTCPP},  $\Pi^{r,y}_{s,t}$ is a  time  consistent convex dynamic procedure on $\hat {\cal W}^r_t$.\\
With the notations of Proposition \ref{propsup}, let $\tilde h=\sup_{\gamma \in \Theta(\Gamma)^0_t}h^{\gamma}$. From equation (\ref{eqetasj2}), it follows that for all $s$ and all $0 \leq r \leq s$, 
\begin{equation}
\sup_{\gamma \in \Theta(\Gamma)^r_t}h^{\gamma}(s,x)= \sup_{\eta \in \Theta(\Gamma)^s_t}h^{\eta}(s,x) 
\label{eqtildeh}
\end{equation}
It follows from Proposition \ref{propsup}  and equation (\ref{eqtildeh}) that $\forall r \leq s \leq t$,
\begin{eqnarray}
\Pi^{r,y}_{s,t}(h(X_t))&=&\sup_{Q^{\gamma}_{r,y} \in {\cal Q}^{\theta}(\Gamma)_{r,y}}  (T^{Q^{\gamma}_{r,y}}_s(h(X_t))- \alpha_{s,t}(Q^{\gamma}_{r,y}))\nonumber\\
&=& \sup_{\eta \in \Theta(\Gamma)^r_t} h^{\eta}(s,X_s)=\tilde h(s,X_s) 
\label{eqtcd8}
\end{eqnarray}
The last assertion follows from Theorem \ref{ThmTCCP}.
\section{Viscosity solution }
\label{viscsol}
The pocedures are the procedures constructed in the previous Section in case of Diffusions with Levy generator or in case of continuous diffusions.
From now on, one assumes furthermore that  $\Gamma$ is a multivalued  Borel mapping convex and closed valued.
\subsection{Viscosity supersolution }
One  assumes that  for all $K $ large enough, $\Gamma_K$ is lower hemicontinuous (cf Definition 16.2 in \cite{AB}, where $\Gamma_K(t,x)=\{(a,b,\mu) \in \Gamma(t,x),\;||a|| \leq K,||b|| \leq K, \int_{\R^n-\{0\}} [||y||^2 1_{||y|| \leq 1}+||y||1_{||y|| >1}]d\mu(y) \leq K\}$ in case of diffusions with Levy generator (and $\Gamma_K(t,x)=\{(a,b) \in \Gamma(t,x),\;||a|| \leq K,||b|| \leq K$ in case of continuous diffusions).
 Let $(t_0,x_0) \in [0,t[ \times \R^n$.
 Let $\phi \in {\cal C}^{1,2}_b([0,t] \times \R^n )$  such that 
\begin{equation}
0=v(t_0,x_0)-\phi(t_0,x_0)=\min(v(t,x)-\phi(t,x))
\label{eqinphi}
\end{equation}
\begin{lemma}
Let $\phi \in {\cal C}^{1,2}_b([0,t] \times \R^n )$. Assume that  $M$ satisfies the hypothesis $M_C$ of Definition \ref{deftetacont2}.  Let 
\begin{equation} 
\tilde K\phi(s,x)(y)=\phi(s,x+y)-\phi(s,x)-\frac{y^*\nabla \phi(s,x)}{1+||y||^2}
\label{eqtilk}
\end{equation}
The function  
$K\phi(s,x)=\int_{\R^n-\{0\}}\tilde K\phi(s,x)(y)M(s,x, dy)$ is a continuous function of $(s,x)$.
\label{lemmakcon}
\end{lemma}
{\bf Proof}
The functions  
$\tilde K\phi(s,x)$ given by (\ref{eqtilk}) are  uniformly bounded. Thus  for all $\epsilon>0$, there is $K>0$ such that for all $\mu \in {\cal M}_C$ (Definition \ref{defMc}), 
\begin{equation}
\forall (s,x),\;\int_{||y|| \geq K}  \tilde K \phi(s,x)(y)\mu(dy)< \epsilon
\label{eq1}
\end{equation}
It follows from the Taylor formula with integral remainder that
$\tilde K \phi(s,x)(y)-\tilde K \phi(s_0,x_0)(y)=y^*(\nabla(\phi)(s,x)-\nabla(\phi)(s_0,x_0))\frac{||y||^2}{1+||y||^2}+  \int_0^1Tr(D^2\phi(s,x+uy)-D^2 \phi(s_0,x_0+uy))yy^*](1-u)du$. 
Let $K_0>0$. Making use of the uniform continuity of $(s,x) \rightarrow \nabla(f)(s,x)$ and $(s,x) \rightarrow D^2\phi(s,x)$ on compact spaces it follows that for $|s-s_0|<\eta$, $||x-x_0||<\eta$, $||x|| \leq K_0$ and $||x_0|| \leq K_0$,  $\int_{||y|| \leq K} |\tilde K\phi(s,x,y)-\tilde K\phi(s_0,x_0,y)| \leq \int_{||y|| \leq K} \epsilon||y||^2\mu(dy)$ for all $\mu \in {\cal M}_C$.
Thus  $\int_{\R^n-\{0\}} |\tilde K\phi(s,x,y)-\tilde K\phi(s_0,x_0,y)|\mu(dy)$ tends to $0$ uniformly in $\mu \in {\cal M}_C$ when   $(s,x)$ tends $(s_0,x_0)$.
On the other hand, making use one more time of the Taylor formula and of the continuity hypothesis (Definition \ref{deftetacont2} 2.), it follows that $\int_{\R^n-\{0\}}\tilde K\phi(s_0,x_0,y)M(s,x,dy)$ is a continuous function of $(s,x)$. 
\hfill $\square $ 

\begin{theorem} The hypothesis are those of Theorem \ref{thmfellerpro}. Assume also that for  $K$ large enough, $\Gamma_K$ is lower hemi continuous. Let $h$ be  continuous   bounded from below.    Let $v=\tilde h$ be a lower semi continuous function such that equations (\ref{eqmeas20}) and (\ref{eqmeas0}) are satisfied. In case of Levy diffusions $v$ is a viscosity supersolution of 
\begin{equation}
\left \{\parbox{12cm}{
\begin{eqnarray} 
-\partial_u v(u,x)-f(u, x,Dv(u,x),D^2v(u,x),\tilde K v(u,x))&=&0 \nonumber \\
v(t,x)& = & h(x)\nonumber 
\end{eqnarray}
}\right.
\label{edpvisc}
\end{equation}
at each point $(t_0,x_0)$ such that  $f(t_0,x_0, D\phi(t_0,x_0), D^2\phi(t_0,x_0), \tilde K \phi(t_0,x_0)))<\infty$ for all $\phi \in {\cal C}_b^{1,2}$
satisfying (\ref{eqinphi}).
Here 
$f(t,x, D \phi(t,x),D^2\phi(t,x),\tilde K \phi(t,x))=
\sup_{(a,b,\mu) \in \Gamma(t,x)}[(b^*D\phi(t,x)+\frac{1}{2}Tr(aD^2\phi(t,x))+\int\tilde K \phi(t,x)(y) \mu(dy)-g(t,x,a,b,\mu)]$, and $\tilde K \phi(t,x))$ is given by equation (\ref{eqtilk})
\label{thmsuper0}\\
In case of continuous diffusions, assume  that  the restriction of $g$  to $\{(u,x,y),\; y \in \Gamma(u,x)\}$ is upper semi-continuous. Then 
 $v$ is a viscosity supersolution of 
\begin{equation}
\left \{\parbox{12cm}{
\begin{eqnarray} 
-\partial_u v(u,x)-f(u, x,Dv(u,x),D^2v(u,x))&=&0 \nonumber \label{eqedp000}\\
v(t,x)& = & h(x)\nonumber 
\end{eqnarray}
}\right.
\end{equation}
at each point $(t_0,x_0)$ such that  $f(t_0,x_0, D\phi(t_0,x_0),D^2\phi(t_0,x_0))<\infty$ for all $\phi \in {\cal C}_b^{1,2}$
satisfying (\ref{eqinphi}).
In this case 
\begin{equation}
f(t,x, z,\gamma)=
\sup_{(a,b) \in \Gamma(t,x)}[(b^*z+\frac{1}{2}Tr(a\gamma)-g(t,x,a,b)]
\label{eqfjo}
\end{equation}
\label{thmsuper}
\end{theorem}
 {\bf Proof} The proof follows the proof of Theorem 5.9 of \cite{JBN-PDE}, replacing Ito's formula by the martingale property.
 From the time consistency for the process $\Pi^{t_0,x_0}_{u,s}$, for all $\delta>0$, 
\begin{eqnarray}
v(t_0,x_0)=\Pi^{t_0,x_0}_{t_0,t}(h(X_t))=\Pi^{t_0,x_0}_{t_0,t_0+\delta}(\Pi^{t_0,x_0}_{t_0 +\delta,t}(h(X_t)))
\end{eqnarray}
From the definition of $v$  and  the inequality $v \geq \phi$, it follows that  \begin{equation}
\Pi^{t_0,x_0}_{t_0+\delta,t}(h(X_t))=v(t_0+\delta, X_{t_0+\delta}) \geq \phi(t_0+\delta, X_{t_0+\delta})
\label{ineq02}
\end{equation}
From the definition of $\Pi^{t_0,x_0}_{t_0,t_0+\delta}$, it follows that for all $\theta \in \Theta(\Gamma)$,
\begin{equation}
v(t_0,x_0) \geq Q^{\theta}_{t_0,x_0}(v(t_0+\delta, X_{t_0+\delta}))-\alpha_{t_0,t_{0+\delta}}( Q^{\theta}_{t_0,x_0})
\label{ineq03}
\end{equation} 
\begin{itemize}
\item We give details in case of diffusions with Levy generator.
Let 
$\tilde K\phi(s,x)$ be given by equation  (\ref{eqtilk}).
The function $\phi$ belongs to  ${\cal C}^{1,2}_b$. Thus it follows from  the martingale property for $Q^{\theta}_{t_0,x_0}$ as stated in Theorem 1.1 of \cite{St}, from the definition of the penalty (Definition \ref{defpenlev}) and the equality $v(t_0,x_0)=\phi(t_0,x_0)$ that  
\begin{eqnarray}
0 \geq &Q^{\theta}_{t_0,x_0}[\int_{t_0} ^{t_0+\delta}\{\partial_u \phi(u,X_u)+\frac{1}{2}Tr(aD^2\phi)(u,X_u)+b^*D\phi(u,X_u)]du\nonumber\\
&+Q^{\theta}_{t_0,x_0}[\int_{t_0} ^{t_0+\delta} [ \int_{\R^n-\{0\}}(\tilde K \phi(u,X_u)(y) M(u,X_u, dy))\}du ]\nonumber\\ &+ Q^{\theta}_{t_0,x_0}[(\int_{t_0} ^{t_0+\delta} -g(u,X_u,a(u,X_u),b(u,X_u),M(u,X_u)du ]
\label{ineq04}
\end{eqnarray} 
$f(t_0,x_0, D\phi(t_0,x_0),D^2\phi(t_0,x_0),\tilde K \phi(t_0,x_0))=
\sup_{(a,b,\mu) \in \Gamma(t_0,x_0)}[b^*D\phi(t_0,x_0)+\frac{1}{2}Tr(aD^2\phi(t_0,x_0))+\int_{\R^n-\{0\}}\tilde K \phi(t_0,x_0)(y) \mu(dy)-g(t_0,x_0,a,b,\mu)]$. Let $(a_0,b_0,\mu_0) \in \Gamma$ such that 
\begin{eqnarray}
&(b_0^*D\phi(t_0,x_0)+\frac{1}{2}Tr(a_0D^2\phi(t_0,x_0))+\int \tilde K\phi(t_0,x_0)d\mu_0-g(t_0,x_0,a_0,b_0,\mu_0)]\nonumber\\
&>f(t_0,x_0, D\phi(t_0,x_0),D^2\phi(t_0,x_0),\tilde K \phi(t_0,x_0))-\epsilon
\label{eqcp2}
\end{eqnarray}
Let $K$ such that $(a_0,b_0,\mu_0) \in \Gamma_K$ and 
$\Gamma_K$ is lower hemicontinuous.  There is then a continuous function $(s,x) \rightarrow a(s,x),b(s,x),M(s,x))=\theta^0(s,x) \in \Gamma_K$ such that $a(s_0,x_0),b(s_0,x_0),M(s_0,x_0)=(a_0,b_0,\mu_0)$. \\
 From the regularity properties of $\phi$ and from Lemma \ref{lemmakcon}, one can write (\ref{ineq04}) as
\begin{equation}
0 \geq Q^{\theta}_{t_0,x_0}[\int_{t_0} ^{t_0+\delta}(\xi(u,X_u)-g^{\theta^0}(u,X_u))du
\label{eqgLevy}
\end{equation}
 where $\xi$ and $g^{\theta^0}$ are  continuous bounded. 
\item In case of continuous diffusions, it follows from equation (\ref{ineq03}) that 
\begin{equation}
0 \geq Q^{\theta}_{t_0,x_0}[\int_{t_0} ^{t_0+\delta}(\xi(u,X_u)-\overline g(u,X_u))du
\label{eqgLevc}
\end{equation}
where $\xi$ is   continuous bounded and $\overline g$ is upper semi-continuous. 
\end{itemize}
For all $\epsilon>0$, there is $\eta>0$,  such that for $t_0 \leq u \leq t \leq t_0+\eta$ and $||x_0-x||<\eta$, 
\begin{eqnarray}
\xi(u, x)-\tilde g(u,x) \geq \xi(t_0, x_0)-\tilde  g(t_0,x_0)-\epsilon 
\label{eqsubv}
\end{eqnarray}
with $\tilde g$ equal $g^{\theta^0}$ or $\overline g$.
In case of Levy generators, let $K_0=||\xi||_{\infty}+||g^{\theta^0}||_{\infty}$.
In case of continuous diffusions, from hypothesis $H_g$, $|\overline g(u,x)|\leq 1+||x||^m$.  Let $K_0=||\xi||_{\infty}+1 +(Q^{\theta_0}_{t_0,x_0}(\sup_{t_0 \leq u \leq \eta}||X_u||^{2m})^{\frac{1}{2}}$.

From Lemma \ref{lemmaLM} in case of Levy generators and from Proposition 2.3 of \cite{JBN-PDE} in case of continuous diffusions,  applied with the probability measure $Q^{\theta_0}_{t_0,x_0}$ there is  $0<\delta < \eta$  such that 
\begin{equation}
Q^{\theta_0}_{t_0,x_0}(A)<\frac{\epsilon}{K_0}\;\; with \;\;A=\{\omega \;|\; \sup_{t_0 \leq u \leq t_0+ \delta}||X_u-x_0||>\eta\;\}
\label{eqcpC}
\end{equation}
Dividing the inequality (\ref{ineq04}) by $\delta$, making use of (\ref{eqcp2}), (\ref{eqsubv}),  (\ref{eqcpC}) and of Cauchy Schwarz inequality, it follows that  
$$0>f(t_0,x_0, D\phi(t_0,x_0,D^2\phi(t_0,x_0)\tilde K \phi(t_0,x_0))-3\epsilon$$.
\hfill $\square $ 
 \begin{remark}
A different  approach to viscosity  solutions of   fully non linear PDE is proposed in \cite{STZ2}. This approach is based on the existence of solutions to BSDE.  In \cite{STZ2} the PDE is $$-\partial_u v(u,x)-f(u, x,v(u,x),Dv(u,x),D^2v(u,x)) $$ with  
$f(t,x, y,z,\gamma)=
\sup_a[\frac{1}{2}Tr(a\gamma)-g(t,x,y,z,a)]$, $f$ is allowed to depend on $y$. However unlike for 
equation (\ref{eqfjo}), the assumptions on function  $g$ (Assumption 5.1 of \cite{STZ2}) are quite restrictive:   the domain of  $g(t,x,.)$ is  independent of $t$ and $x$, the function $g$ is   uniformly continuous  in $t$, Lipschitz in $z$.
\end{remark}
%====================================================================================
\subsection{Viscosity subsolution}
From now on we add a new hypothesis $H'$. In case of continuous diffusions we assume that the multivalued Borel mapping $\Gamma$ has linear growth, i.e. There is a constant $K>0$ such that $ \forall (a,b) \in \Gamma(s,x),\; ||a||,||b|| \leq K(1+||x||)$.  As in \cite{JBN-PDE} Lemma 5.11, it follows from \cite{Kry}  II 5 Corollary 10 that 
for all $q\geq 1$,  $C$,  and $t>0$, there is a   constant $K_1$ such that for all $y$ such that $||y|| \leq C$  and all bounded $\gamma \in \Theta(\Gamma)^r_t$, 
\begin{equation}
E_{Q^{\gamma}_{r,y}}(\sup_{s \leq u \leq t }(||X_t-X_{u}||^{2q}) \leq K_1 (t-s)^q
\label{eqKry01}
\end{equation}
and 
\begin{equation}
\forall \epsilon>0,\forall \eta>0,\exists h>0 \forall r>0, \forall y, ||y||  \leq C,\; \forall \gamma\in \Theta(\Gamma)^r_t, \; Q^{\gamma}_{r,y}[\sup_{r \leq v \leq r+h}||X_v-y||>\eta)< \epsilon
\label{equnif}
\end{equation}

Denote $v^*$ the upper semi continuous envelope of $v$, $$v^*(s,x)=\limsup_{(s',x')\rightarrow (s,x)}    v(s',x')$$
Let $\phi \in {\cal C}^{1,2}_{b}$ such that 
\begin{equation}
0=v^*(t_0,x_0)-\phi(t_0,x_0)=\sup_{(s,x)}(v^*(s,x)-\phi(s,x))
\label{equsc}
\end{equation}
\begin{theorem} Let  $h$ be  continuous bounded from below. Assume that 
\begin{enumerate}
\item  In case of continuous diffusions, $g$  satisfies hypothesis $H_g$ and $\Gamma$ is closed convex with linear growth.
\item  In case of diffusions with Levy generator,  $\Gamma=\Gamma_K$ for some $K>0$. 
\end{enumerate}
Assume that $f$ is upper semi continuous. 
Let $v=\tilde h$ be the lower semi continuous function  as in Theorem \ref{thmfellerpro}. Then $v^*$ is a viscosity subsolution of (\ref{eqedp00}) in case of continuous diffusions and of (\ref{edpvisc}) in case of diffusions with Levy generator
\label{thmsub}
\end{theorem}
$\tilde K \phi(u,x)$ belongs to the vector space ${\cal C}(\R^n)$ of continuous functions on $\R^n$. Let ${\cal M}$ be the vector space of measures generated by ${\cal M}_C$. The topology on $E$ is $\sigma(E,{\cal M})$.\\ 
{\bf Proof} It follows from the first part of the proof of Lemma \ref{lemmakcon}  that $\tilde K \phi$ is continuous from $[0,T] \times \R^n$ to $E$. The functions $\phi_u$,$D\phi$, $D^2\phi$  are continuous and $f$ is upper semi-continuous thus  for all $n>0$ there is $\eta_n>0$ such that for $t_0 \leq u \leq t_0+\eta_n$, and $||x-x_0||<\eta_n$, $\phi_u(u,x)+f(u,x, D\phi(u,x), D^2\phi(u,x),\tilde K \phi(u,x) \leq$\\ 
$\phi_u(t_0,x_0)+f(t_0,x_0, D\phi(t_0,x_0), D^2\phi(t_0,x_0),\tilde K \phi(t_0,x_0)+\frac{1}{n}$.
From equation (\ref{equnif}) in case of continuous diffusions and inequality (\ref{eqLM})  in case of diffusions with Levy generator,
there is $\alpha_n>0$ such that for all $|u-t_0|\leq 1$ and $||y-x_0|| \leq 1$ for all $\gamma \in \Theta(\Gamma)^u_t$, 
$${Q^{\gamma}_{u,y}}(A_n)<\frac{1}{n}\;\; \text{with} \;\;A_n=\{\omega \;|\; \sup_{u \leq u' \leq u+ \alpha_n}||X_{u'}-y||>\frac{\eta_n}{2}\;\}$$
Let $\delta_n=inf(\eta_n,\alpha_n)$.
For all $n>0$ choose $(t_n,x_n)$ such that  $|t_n-t_0|<\frac{\delta_n}{2}$, $||x_n-x_0||<\frac{\eta_n}{2}$ and $\phi(t_n,x_n)-  v(t_n,x_n) \leq \frac{\delta_n}{n}$. 
For all $n>0$, there is a  process $\gamma_n \in \Theta(\Gamma)^{t_n}_t$, 
such that 
\begin{equation}
v(t_n,x_n) \leq E_{Q^{\gamma_n}_{t_n,x_n}}(h(X_t))-\alpha_{t_n,t}(Q^{\gamma_n}_{t_n,x_n})+\frac{1}{n}.
\label{eq24}
\end{equation}
From the cocycle condition for the penalty associated to the probability measure $Q^{\gamma_n}_{t_n,x_n}$, the definition of  $\Pi^{t_n,x_n}_{t_n+\delta_n,t}(h(X_t)$ and $\Pi^{t_n,x_n}_{t_n+\delta_n,t}(h(X_t)=v(t_n + \delta_n,X_{t_n + \delta_n}) \leq \phi(t_n + \delta_n,X_{t_n + \delta_n})$
\begin{equation}
v(t_n,x_n) \leq E_{Q^{\gamma_n}_{t_n,x_n}}[\phi(t_n+\delta_n,X_{t_n+\delta_n})-\int_{t_n}^{t_n+\delta_n}g (u,X_u,\mu_n(u,\omega)du]+\frac{1}{n}
\label{eq25}
\end{equation}
We apply  the martingale property for $Q^{\gamma_n}_{t_n,x_n}$ and we make use of  the definition of $f$.
We can then finish the proof as the proof of the viscosity supersolution. 
 
\subsection{Viscosity solution and uniqueness}
\label{subsecU}
The following Theorem results from Theorems \ref{thmsuper} and \ref{thmsub}
\begin{theorem}
Assume that all the previous hypothesis are satisfied. Assume furthermore that $h $ is   bounded from below and $\alpha$ H\"older-continuous for some $\alpha>0$ in case of continuous diffusions  and $h$ is Lipschitz in case of diffusions with Levy generator. Let $v=\tilde h$ be a lower semi continuous function such that equations (\ref{eqmeas20}) and (\ref{eqmeas0}) are satisfied.  Assume that the PDE (\ref{edpvisc}) satisfies the comparison principle for functions bounded on compact spaces. Then the function $v$ is continuous. It is the unique viscosity solution of  (\ref{edpvisc}). 
\end{theorem}
{\bf Proof}
As in \cite{JBN-PDE} we prove first that $v$ is continuous at $(t,x)$ for all $x$. Making  use of the equation (\ref{eqKry01}) in case of continuous diffusions and of Lemma  \ref{lemmaLM} for diffusions with levy generator, we prove as in \cite{JBN-PDE} that $-v(s,x)-h(x)$ tends to $0$  uniformly in $x$ when $s$ tends to $t$. The result follows then from Thm \ref{thmsuper} and \ref{thmsub} and from Proposition 5.5 of \cite{JBN-PDE}.
For comparison results for non linear second order PDE we refer to \cite{GIL,FSo,Da}.\\


\begin{thebibliography}{777}

\bibitem{AB} Aliprantis C., Border K., Infinite dimensional analysis, 2nde edition, Springer Berlin, 1999.

\bibitem{Bil} Billingsley P., Convergence of Probability Measures, Wiley Series in Probability and Mathematical Statistics, Wiley and  Sons, New York, 1999.

\bibitem{BN03}
Bion-Nadal J., Dynamic risk measures: Time consistency and risk measures
from BMO martingales, Finance and Stochastics 12  (2008)

\bibitem{BN04}
Bion-Nadal J. , Time Consistent Dynamic Risk Processes,  Stochastic Processes and their Apllications,   119 (2009), 633-654.

\bibitem{JBN-PDE}
Bion-Nadal J., Time consistent convex Feller processes and non linear second order partial differential equations, preprint arXiv:1207.1742

\bibitem{BNK1} Bion-Nadal J. and Kervarec M., Risk Mesuring under Model Uncertainty, The annals of applied probability, 22 (2012), 213-238

\bibitem{BNK2}
Bion-Nadal J. and Kervarec M., Dynamic risk measuring under model uncertainty: taking advantage of the hidden probability measure, preprint arXiv:1012.5850

\bibitem{CDK}
Cheredito P., Delbaen F. and Kupper M., Dynamic monetary risk measures for bounded discrete time processes, Electronic Journal of Probability,  \textbf{11} (2006)  57-106 

\bibitem{Da}
Da Lio F. and Ley O. Convex Hamiltonian-Jacobi equations under superlinear growth conditions on data, Appl. Math. Optim. \textbf{63} (3) (2011) 309-339



\bibitem{D}
Delbaen F., The structure of $m$-stable sets and in particular of the set of risk neutral measures, In Memoriam Paul-Andr\'e Meyer, in: Lecture Notes in Mathematics, \textbf{1874} (2006)  215-258 

\bibitem{DPR}
Delbaen F., Peng S. and Rosazza Gianin E.: Representation of the penalty term of dynamic concave utilities, Finance and Stochastics, 14, pp 449-472 (2010)


\bibitem{DHP} Denis L., Hu M. and Peng S., Function spaces and capacity related to a Sublinear Expectatio: application to G-Brownian Motion Pathes, Potential Analysis \textbf{34}(2) pp 139-161 (2011) %% preprint arXiv:0802.1240

\bibitem {EK2}
El Karoui N, Lepeltier J.P. and Millet A., A probabilistic approach of the reduite, probability and mathematical Statistics \textbf{13} (1) (1992), pp 97-121.

\bibitem{FSo} Fleming W.H. and Soner H. M. Controlled Markov Processes and Viscosity Solutions, 2nd edition, Springer (2006)

\bibitem{FS}
F\"ollmer H. and Schied A., Stochastic Finance: An introduction in discrete time, second edition, de Gruyter studies in mathematics 27 (2004).


\bibitem{FR}
Frittelli M., Rosazza Gianin E., Putting Order in Risk Measures. Journal of Banking Finance 26, (2002) pp. 1473-1486 


\bibitem{GIL}
Grandall M., Ishii H. and Lions P-L: User's guide to viscosity solutions of second order partial differential equations. Bulletin of the american mathematical society \textbf{27}, Number 1, pp  1-67 (1992)

\bibitem{HP}
Hu M. and Peng S., G-L\'evy processes under sublinear expectations,    arXiv:0911.3533

\bibitem{KPZ}
Kazi-Tani N., Possamai D. and Zhou C., Second order BSDE with jumps. arXiv:1208.0763

\bibitem{KS} 
Kl\"oppel S. and Schweizer M.,  Dynamic indifference valuation via convex risk measures, %%  preprint ETH, 2005.
Mathematical Finance, Volume 17, Issue 4, pages 599–627, October 2007

\bibitem{Kry} 
Krylov N., Controlled Diffusion Processes, Springer (1980)

\bibitem {LM}
Lepeltier J.P. and Marchal B., Probl\`eme des martingales et \'equations diff\'erentielles stochastiques associ\'ees \`a un op\'erateur int\'egro diff\'erentiel, Annales de l'IHP, section B, \textbf{22} (1) pp 43-103, (1976)

\bibitem{N1} Nutz M., Random G-expectations, to appear in Annals of Applied Probability

\bibitem{N2} Nutz M., A quasi-sure approach to the control of non-Markovian stochastic differential equations, Electronic Journal of Probability, \textbf{17}, No. 23, pp. 1-23, (2012)

\bibitem{NS}
Nutz M. and Soner H.M., Superhedging and dynamic risk measures under volatility uncertainty, SIAM Journal on Control and Optimization, \textbf{ 50}, No. 4, pp. 2065-2089, (2012)

\bibitem{NVH} Nutz M. and van Handel R., Constructing sublinear expectations on path space, preprint arXiv:1205.2415

\bibitem {P1}
Peng S., G-expectation, G-Brownian motion and related stochastic calculus of Ito Type, in  Stochastic Analysis and A, vol 2 of Abel Symp., (2007), pp 541-567, Springer, Berlin.

\bibitem {P2}
Peng S., Multi-dimensional  G-Brownian motion and related stochastic calculus underv G-expectation,   Stochastic Processes and their Apllications,   118 (2008), pp 2223-2253. 



\bibitem{RY}
Revuz D. and Yor M., Continuous martingales and Brownian motion, Springer.

\bibitem{STZ1} Soner H.M. Touzi N. and Zhang J., Dual formulation of second oredr target problems, to appear in annals of applied probability

\bibitem{STZ2} Soner H.M. Touzi N. and Zhang J., Wellposedness of second order backward SDE, to appear in Probability Theory and Related Fields 

\bibitem{St} 
Strook D. Diffusion processes asociated with Levy generators, Z. Wahrscheinlichkeitstheorie verw. Gebiete \textbf{32} (1975) 209-244

\bibitem{SV1} 
Stroock and Varadhan, Diffusion processes with continuous coefficients I, Communications on Pure and Applied Mathematics, \textbf{22}, pp 345-400 (1969).

\bibitem{SV2} 
Stroock and Varadhan, Diffusion processes with continuous coefficients II, Communications on Pure and Applied Mathematics, \textbf{22}, pp 479-530 (1969).


\bibitem{SV3} 
Stroock and Varadhan, Diffusion processes with boundary conditions. Comm. Pure Apll. Math \textbf{24}, (1971) 147-225.

\bibitem{SV4}
Stroock D.W. and Varadhan S.R.S.: Multidimensional diffusion processes, Springer, 1979. 

\end{thebibliography}
\end{document}